\documentclass[11pt,twoside]{amsart}

\usepackage[utf8]{inputenc}
\usepackage[T1]{fontenc}
\usepackage{amsmath,amssymb,amsthm,mathtools}
\usepackage[margin=1in]{geometry}
\usepackage{tikz-cd}
\usepackage[colorlinks=true,linkcolor=blue,citecolor=blue]{hyperref}
\usepackage{cleveref}

\theoremstyle{plain}
\newtheorem{theorem}{Theorem}[section]
\newtheorem{proposition}[theorem]{Proposition}
\newtheorem{lemma}[theorem]{Lemma}
\newtheorem{corollary}[theorem]{Corollary}
\newtheorem{claim}[theorem]{Claim}

\theoremstyle{definition}
\newtheorem{definition}[theorem]{Definition}
\newtheorem{remark}[theorem]{Remark}

\newtheorem{setup}[theorem]{Set-up}

\newcommand{\PP}{\mathbb{P}}
\newcommand{\PM}{\mathbb{P}^M}
\newcommand{\OO}{\mathcal{O}}
\newcommand{\sHom}{\mathcal{H}om}
\newcommand{\Hom}{\operatorname{Hom}}
\newcommand{\Ext}{\operatorname{Ext}}
\newcommand{\coker}{\operatorname{coker}}
\renewcommand{\Im}{\operatorname{Im}}
\newcommand{\id}{\mathrm{id}}
\newcommand{\Yt}{\widetilde{Y}}
\newcommand{\Ht}{\widetilde{H}}
\newcommand{\Ct}{\widetilde{C}}
\newcommand{\yt}{\tilde{y}}
\newcommand{\NY}{N_{Y}}
\newcommand{\NYt}{N_{\Yt}}
\newcommand{\Hc}{\mathcal{H}}
\newcommand{\TP}{T_{\PM}|_Y}

\newcommand{\SC}{\mathcal{S}_C}
\newcommand{\RC}{\mathcal{R}_C}
\newcommand{\GH}{G_H}
\newcommand{\Dd}{D}
\newcommand{\cY}{\mathcal{Y}}
\newcommand{\cYt}{\widetilde{\mathcal{Y}}}
\newcommand{\cCt}{\widetilde{\mathcal{C}}}
\newcommand{\cN}{\mathcal{N}}
\newcommand{\pr}{\mathrm{pr}}
\makeatletter
\providecommand{\@wraptoccontribs}[2]{}
\def\@setauthors{%
  \begingroup
  \def\thanks{\protect\thanks@warning}%
  \trivlist
  \centering\footnotesize \@topsep30\p@\relax
  \advance\@topsep by -\baselineskip
  \item\relax
  \author@andify\authors
  \def\\{\protect\linebreak}%
  \MakeUppercase{\authors}%
  \ifx\@empty\contribs
  \else
    \\ \@setcontribs
    \@closetoccontribs
  \fi
  \endtrivlist
  \endgroup
}
\makeatother
\title[Extendability of general prime $K3$ surfaces without Gaussian maps ]{Extendability of general prime $K3$ surfaces without Gaussian maps and classification of prime Fano threefolds}
\author{Jayan Mukherjee}
\address{Department of Mathematics, Oklahoma State University}
\email{jayan.mukherjee@okstate.edu}
\contrib[(with an appendix written jointly with]{Anand Deopurkar)}
\address{Mathematical Sciences Institute, Australian National University}
\email{anand.deopurkar@anu.edu.au}

\begin{document}

\begin{abstract}
In \cite{BM24} an approach to the question of extendability of projective varieties via degeneration to ribbons was introduced, and in \cite{BM25} it was used to give a new proof of optimal results on the extendability of general non-prime $K3$ surfaces and the classification of non-prime Fano threefolds. In this article we treat the prime case. We give a new proof of the optimal results on the extendability of general prime $K3$ surfaces of every genus, of the classification of prime Fano threefolds of Picard rank one due to Fano, Iskovskikh and Mukai, and of the irreducibility of their Hilbert schemes. This involves degenerations of $K3$ surfaces to ribbons on embedded rational normal scrolls, called $K3$ carpets. We compute exactly the cohomology of the twisted normal bundle of these $K3$ carpets, instead of computing coranks of Gaussian maps of the canonical curve sections as in \cite{CLM93}. The key challenge is to analyze the coboundary map of the sequence relating the normal bundle of the $K3$ carpet to the normal bundle of the scroll. We give a deformation theoretic interpretation of this map. Given such a $K3$ carpet and its hyperplane section, which is a canonical ribbon, we show that the coboundary map corresponds to the obstruction to lifting a first order deformation of the underlying scroll fixing the underlying rational normal curve to a first order deformation of the $K3$ carpet fixing the canonical ribbon. To compute the obstruction, we use a result on extensions of non-split ribbons on any curve to a unique $K$-trivial carpet in a canonical way. The genus zero case of this result was stated in \cite{BE95}. %This result in fact independently proves the birationality  of the flag moduli space $(X,C)$, where $X$ is a $K3$ surface and $C$ is a canonical curve to its image under the second projection, as predicted by Mukai's program of constructing the $K3$ surface from the canonical curve. 
%Together with \cite{BM24} and \cite{BM25}, this completes the proof of the extendability results for general $K3$ surfaces and of the classification of Fano threefolds of Picard rank one by means of ribbon degenerations of Hirzebruch surfaces, with no recourse to Gaussian maps.
\end{abstract}

\maketitle

\section{Introduction}

\begin{definition}\label{extendability}
Let $X \subset \mathbb{P}^M$ be an irreducible nondegenerate variety of codimension at least $1$. Let $k \geq 1$ be an integer. We say that $X$ is $k-$ extendable if there exists a variety $W \subset \mathbb{P}^{M+k}$
different from a cone, with dim $W =$ dim $X + k$ and having $X$ as a section by a $M-$ dimensional linear space such that $W$ is smooth along $X = W \cap \mathbb{P}^M$. We say that $X$ is precisely $k-$ extendable if it is $k-$ extendable but not $(k + 1)-$ extendable. The variety $W$ is called a $k$-extension of $X$. We say that $X$ is extendable if it is $1-$ extendable.
\end{definition}

Extendability of a projective variety is a natural and fundamental question in projective geometry. This classical question has been the topic of intense research for decades and has revealed deep connections between the geometry of an embedding, Gaussian map of curve sections, deformations of cones over the hyperplane sections etc. We refer to \cite{Lop23} for an excellent recent survey of this topic.

\smallskip

This article revisits the question of extendability of an embedded prime $K3$ surface and its applications to the boundedness, classification and Hilbert schemes of prime Fano threefolds of Picard rank one. In \cite{CLM93}, \cite{CLM94} and \cite{CLM98}, Ciliberto, Lopez and Miranda analyzed the extendability of general $K3$ surfaces by computing the coranks of Wahl maps of the canonical curve sections (see also \cite{ABS17}, \cite{CD22}, \cite{CD24}, \cite{CDGK20}, \cite{CDS20}, \cite{K01}, \cite{K20}, \cite{KLM11}, \cite{Tot20} for related and more recent work). 
%In the case of prime $K3$ surfaces, which is the case treated in this article, the authors of \cite{CLM93} degenerate the $K3$ surface to a union of two rational normal scrolls, and further to a union of planes whose hyperplane sections are graph curves, and compute in this way the corank of the Wahl map of a general canonical curve section of a general prime $K3$ surface of genus $g$. It turns out to be one if and only if $g = 11$ or $g \geq 13$. The classification of prime Fano threefolds of Picard rank one, due to Fano, Iskovskikh and Mukai (see \cite{Isk77}, \cite{Isk78}, \cite{Muk88}), is then recovered by a careful analysis of the Hilbert scheme at the points representing cones over canonical curves (see \cite{CLM93}, \cite{CLM98}).

\smallskip

The questions on extendability has almost always involved the Gaussian map, either on the surface or going to the hyperplane curve section. As noted in \cite{Lop23} ``the emerging philosophy is that if one can control the corank of the Gaussian maps of the curve section of X, then one can give information on the extendability of X.'' 
%It is further noted in \cite{Lop23} that there are only two other general results (\cite{KLM11}, \cite{CDGK20}) that allows studying the extendability of surfaces, without passing to the hyperplane section.

\smallskip

In \cite{BM24} an approach to the question of extendability via degeneration to ribbons was introduced. It avoids Gaussian maps altogether, either on the surface or on the hyperplane curve section, and provides a general method to study extendability, by directly computing $\alpha(X) = h^0(N_{X/\mathbb{P}^M}(-1))-M-1$ and $h^0(N_{X/\mathbb{P}^M}(-2))$ on a ribbon degeneration of $X$. Once $\alpha(X)$ and $h^0(N_{X/\mathbb{P}^M}(-2))$ are computed, we conclude using the following theorem due to Zak-L\'vovsky.

\begin{theorem}(\cite{Lvo92}, \cite{Zak91})\label{ZL}
Let $X \subset \mathbb{P}^{M}$ be a smooth irreducible nondegenerate variety of codimension at least $1$ and suppose $X$ is not a quadric. If $\alpha(X) \leq 0$, then X is not extendable. Further, given an integer $k \geq 2$, suppose that either:
\begin{itemize}
    \item[(1)] $\alpha(X) < M$ or
    \item[(2)] $H^0(N_{X/\mathbb{P}^M}(-2)) = 0$,
\end{itemize}
If $\alpha(X) \leq k-1$, then X is not $k-$ extendable.
\end{theorem}

In \cite{BM25} this method was used to give a new proof of the optimal results of \cite{CLM98} on the extendability of general non-prime $K3$ surfaces and of the classification of non-prime Fano threefolds of Picard rank one. 
%the prime case could be handled there only for infinitely many, but sporadic, values of $g$. 
In this article we treat the prime case and thereby complete this program, i.e., we prove the optimal results on the extendability of general prime $K3$ surfaces, and the classification of prime Fano threefolds of Picard rank one.

By semicontinuity it is enough to give sharp upper bounds to $\alpha(X)$ by computing it over an embedded degeneration of $X$. Let $Y = \mathbb{F}_e$, $e \in \{0,1\}$, be a Hirzebruch surface embedded in $\mathbb{P}^g$ by the complete linear series of $H = C_0+bf$, $b \geq e+2$, as a rational normal scroll of degree $g-1$, where $g = 2b-e+1$. Every genus $g \geq 6$ arises in this way. Let $\widetilde{Y} \subset \mathbb{P}^g$ be the $K3$ carpet on $Y$ (see Definition~\ref{ribbon}), that is, the unique double structure on $Y$ with trivial dualizing sheaf and $h^1(\mathcal{O}_{\widetilde{Y}}) = 0$ \cite{GP97}, and let $\widetilde{H} = \mathcal{O}_{\widetilde{Y}}(1)$. By \cite{GP97} (see also \cite{BMR21}, \cite{Gon06}, \cite{GGP10}, \cite{GGP13}), $\widetilde{Y}$ is smoothable in $\mathbb{P}^g$ to prime $K3$ surfaces of genus $g$ (see Theorem~\ref{degeneration into $K3$ carpets}), so that $\alpha(X) \leq \alpha(\widetilde{Y}) := h^0(N_{\widetilde{Y}/\mathbb{P}^g}(-\widetilde{H}))-g-1$ for a general prime $K3$ surface $X$ of genus $g$, and similarly for $h^0(N_{X/\mathbb{P}^g}(-2))$ (see \eqref{eq:semicont}). It is worth noting that by \cite{BGM24}, such a $\widetilde Y$ lies in the closure of a union of two scrolls meeting transversally along an anti-canonical elliptic curve, a degeneration that was used in \cite{CLM93}.  As in \cite{BM24} and \cite{BM25}, the computability of $\alpha(\widetilde{Y})$ is facilitated by the fact that $\widetilde{Y} \subset \mathbb{P}^g$ is an embedded ribbon which is a local complete intersection scheme whose reduced part is the scroll $Y \subset \mathbb{P}^g$. Three exact sequences (see Remark~\ref{rem:sequences}) reduce the computation of the cohomology of twists of the normal bundle $N_{\widetilde{Y}/\mathbb{P}^g}$ to the computation of the cohomology of twists of $N_{Y/\mathbb{P}^g}$ and to the analysis of the related coboundary maps. The new difficulty in the prime case is that one of these coboundary maps cannot be controlled by cohomological means alone. We give a deformation theoretic interpretation of this map, which we explain in \S\ref{subsec:intro-ingredient}; this is the key new ingredient of this article.

\smallskip

We now state our results. Let $\mathcal{H}_g$ denote the Hilbert scheme of prime $K3$ surfaces of genus $g$ (see Definition~\ref{Hilbert schemes of K3 surfaces} and \S\ref{subsec:fano-setup}), i.e., of smooth $K3$ surfaces $X \subset \mathbb{P}^g$ embedded by the complete linear series of a primitive very ample line bundle $B$ with $B^2 = 2g-2$. It is irreducible, and a point of $\mathcal{H}_g$ is called general if it lies in a suitable dense open subset. Our main theorem computes $\alpha(\widetilde{Y})$ exactly, for $K3$ carpets on rational normal scrolls.

\begin{theorem}(see Theorem~\ref{thm:alpha}, Proposition~\ref{prop:g6}, Proposition~\ref{prop:N2})\label{thm:carpet intro}
Let $Y = \mathbb{F}_e \subset \mathbb{P}^M$, $e \geq 0$, be embedded by the complete linear series of $H = C_0+bf$ with $b \geq e+2$, so that $M = 2b-e+1$, and let $\widetilde{Y} \subset \mathbb{P}^M$ be the $K3$ carpet on $Y$.
\begin{enumerate}
    \item If $b \geq e+3$ or $(e,b) = (1,3)$, then
    \[
    \alpha(\widetilde{Y}) = h^0(\mathcal{O}_Y(-H-2K_Y)) = h^0\big(\mathcal{O}_Y(3C_0+(2e+4-b)f)\big).
    \]
    In particular $\alpha(\widetilde{Y}) = 0$ if and only if $b \geq 2e+5$.
    \item $h^0(N_{\widetilde{Y}/\mathbb{P}^M}(-2\widetilde{H})) \leq h^0(\mathcal{O}_Y(-2H-2K_Y))$. In particular $H^0(N_{\widetilde{Y}/\mathbb{P}^M}(-2\widetilde{H})) = 0$ if $b \geq e+3$, and $h^0(N_{\widetilde{Y}/\mathbb{P}^M}(-2\widetilde{H})) \leq 1$ if $(e,b) = (1,3)$.
\end{enumerate}
\end{theorem}

For general prime $K3$ surfaces we obtain:

\begin{theorem}(see Theorem~\ref{thm:alpha-lowgenus}, Corollary~\ref{cor:cone}, Corollary~\ref{cor:equality})\label{thm:K3 intro}
Let $X \subset \mathbb{P}^g$ be a general prime $K3$ surface of genus $g \geq 6$.
\begin{enumerate}
    \item If $g = 12$, then $\alpha(X) = 1$.
    \item If $g = 10$, then $\alpha(X) = 3$.
    \item If $g = 9$, then $\alpha(X) = 4$.
    \item If $g = 8$, then $\alpha(X) = 6$.
    \item If $g = 7$, then $\alpha(X) = 8$.
    \item If $g = 6$, then $\alpha(X) = 10$.
    \item If $g = 11$ or $g \geq 13$, then $\alpha(X) = 0$ and consequently $X \subset \mathbb{P}^g$ is not extendable.
    \item $H^0(N_{X/\mathbb{P}^g}(-2)) = 0$ if $g \geq 7$, and $h^0(N_{X/\mathbb{P}^g}(-2)) = 1$ if $g = 6$.
\end{enumerate}
\end{theorem}

%The bounds in $(1)$--$(6)$ and $(8)$ are sharp. The proof of Theorem~\ref{thm:fano intro} shows that equality holds in each of them, since the cone over $X$ lies on an irreducible component of the closure of $\mathcal{V}_g$ whose dimension is the upper bound for the dimension of its tangent space given by Theorem~\ref{thm:K3 intro} and Corollary~\ref{cor:cone intro} below. 
Theorem~\ref{thm:K3 intro} recovers the results of \cite{CLM93} on general prime $K3$ surfaces of every genus $g \geq 6$ (for $g \leq 5$ a general prime $K3$ surface is a complete intersection and the results are classical), while the non-prime case was treated in \cite{BM25}.

\smallskip

One of the key technical results proved in \cite{CLM93} and \cite{CLM98}, that crucially depends on the corank of the Wahl map, concerns the dimension of the component of the Hilbert scheme at points representing cones over canonical curves. Using Theorem~\ref{thm:K3 intro} we prove this result directly for the cone over a $K3$ surface.

\begin{corollary}(see Corollary~\ref{cor:cone})\label{cor:cone intro}
Let $X \subset \mathbb{P}^g$ be a general prime $K3$ surface of genus $g \geq 6$ and let $C(X) \subset \mathbb{P}^{g+1}$ be the cone over $X$. Then the tangent space to the Hilbert scheme of $\mathbb{P}^{g+1}$ at $[C(X)]$ has dimension
\[
\dim T_{[C(X)]} = 18+(g+1)^2+(g+1)+\alpha(X)+h^0(N_{X/\mathbb{P}^g}(-2)).
\]
\end{corollary}

The following theorem, which is a consequence of our results and a general principle outlined in \cite{CLM93} and \cite{CLM98}, gives a new proof of the classification of prime Fano threefolds of Picard rank one due to Fano, Iskovskikh and Mukai, and of the irreducibility of their Hilbert schemes.

\begin{theorem}(see Theorem~\ref{thm:fano})\label{thm:fano intro}
Let $g \geq 6$ and let $\mathcal{V}_g$ denote the locus, in the Hilbert scheme of $\mathbb{P}^{g+1}$, of smooth Fano threefolds $V \subset \mathbb{P}^{g+1}$ of Picard rank one, index one and genus $g$, embedded by the complete linear series $|-K_V|$ (see Definition~\ref{def:Vg}).
\begin{enumerate}
    \item $\mathcal{V}_g = \varnothing$ if $g = 11$ or $g \geq 13$.
    \item For $6 \leq g \leq 10$ and $g = 12$, $\mathcal{V}_g$ is irreducible, and the families of Fano and Iskovskikh form an open dense subset of $\mathcal{V}_g$.
\end{enumerate}
\end{theorem}

%Once again, we point out that the above results were proven using very different methods in \cite{CLM93} and \cite{CLM98}. The relation of our degeneration to the one used there is as follows. In \cite{CLM93}, to deal with the prime case, the authors degenerated $K3$ surfaces into a union of two rational normal scrolls meeting along an anticanonical elliptic curve. These were further degenerated to a union of planes whose hyperplane sections are graph curves with corank one Gaussian map. The $K3$ carpet $\widetilde{Y}$ lies in the closure of the locus of such unions of two scrolls (see \cite[Theorem $4.4$]{GP97}). The techniques introduced in \cite{BM24}, \cite{BM25} and in the present article show that degenerating the $K3$ surfaces to the $K3$ carpet suffices to prove the results, requiring no further degeneration to planes, no passage to the hyperplane curve section and no Gaussian maps.

%% ---- optional: uncomment if a corollary on the corank of the Wahl map (as in [BM25]) is added to Section 4 ----
%Finally, exactly as in \cite{BM25}, Theorem~\ref{thm:K3 intro} $(7)$ together with \cite[$2.2$]{CDS20} (see also \cite{ABS17}) and \cite[Lemma $8.3$]{K01} gives back the corank one theorem of \cite{CLM93}: for any smooth hyperplane section $C$ of a general prime $K3$ surface of genus $g = 11$ or $g \geq 13$, one has $\alpha(C) = \operatorname{cork}(\Phi_{\omega_C}) = 1$.

\subsection{The key new ingredient in the proof}\label{subsec:intro-ingredient}

We keep the notation of Theorem~\ref{thm:carpet intro} and assume $b \geq e+3$; the case $(e,b) = (1,3)$ requires a small additional argument (see Proposition~\ref{prop:g6}). We write $\NY := N_{Y/\PM}$ and $\NYt := N_{\Yt/\PM}$. Let $H \subset \PM$ be a hyperplane meeting $Y$ transversally along a smooth curve $C := Y \cap H$ which is a rational normal curve of degree $g-1$ in $H \cong \PP^{g-1}$ (note that $\OO_Y(H) = \OO_Y(C)$). Let $\Ct := \Yt \cap H$, which is a canonical ribbon on $C$ of arithmetic genus $g$ (see \cite{BE95}, \cite{GP97}). Put $N_C := N_{C/H}$ and $L_C := K_Y|_C \cong \OO_{\PP^1}(-g-1)$, the conormal bundle of $\Ct$.

\smallskip

The cohomological procedure of \cite{BM24}, \cite{BM25}, applied to the three exact sequences \eqref{eq:3.1}--\eqref{eq:3.3}, the normal sequence of $Y \subset \PM$ and the Euler sequence, gives (see \eqref{eq:alpha-F})
\[
\alpha(\Yt) = h^0(\OO_Y(-H-2K_Y))+\dim\coker F,
\]
where $F\colon H^0(\NYt\otimes\OO_Y(-H)) \to H^1(\NYt\otimes K_Y(-H))$ is the coboundary map of the exact sequence $0 \to \NYt\otimes K_Y \to \NYt \to \NYt\otimes\OO_Y \to 0$, relating the normal bundle of the $K3$ carpet to its restriction to the scroll, twisted by $\OO_{\Yt}(-\Ht)$. 
%In the non-prime case treated in \cite{BM25} the coboundary maps that appear in this procedure can be controlled by cohomological means (vanishing theorems and pushing forward to $\PP^1$). For the scroll $Y$ this is no longer the case: the target of $F$ is nonzero, and 
Now Theorem~\ref{thm:carpet intro}, $(1)$ is equivalent to the surjectivity of $F$. The surjectivity of $F$ is proved in the following steps.

\begin{enumerate}
    \item[(1)] \emph{Reduction to a map $a$ between cohomology groups of twists of $\NY$} (see the proof of Theorem~\ref{thm:alpha}, Lemma~\ref{lem:hfg} and Proposition~\ref{surjectivity of a}). Using the sequences \eqref{eq:3.2} and \eqref{eq:3.3}, the group $H^0(\NY(-H))$ injects into the source of $F$ and $H^1(\NY\otimes K_Y(-H))$ injects into its target and it is shown that
    %A \v{C}ech cocycle computation (Lemma~\ref{lem:hfg}) shows that $F$ induces a commutative diagram with exact rows \eqref{diag:a} whose right-hand vertical map is $2\cdot\id$; by the snake lemma, 
    $F$ is surjective if and only if the induced map
    \[
    a\colon H^0(\NY(-H)) \longrightarrow H^1(\NY\otimes K_Y(-H))
    \]
    is surjective. Note that $H^0(\NY(-H)) = H^0(\NY(-C))$ is the tangent space $T_{[Y]}\SC$ at $[Y]$ to the space  $\SC$ of smooth rational normal scrolls $Y' \subset \PM$ with $Y' \cap H = C$, which is smooth of dimension $2(g-1)$ at $[Y]$ (see Lemma~\ref{lem:SC}). $\mathcal{S}_C$ is the fiber over $C \subset H$ of the flag Hilbert scheme parametrizing pairs $(Y',C')$, where $Y'$ is a rational normal scroll and $C'$ is a hyperplane section of it, which is a rational normal curve. 
    %A tangent vector is a first order embedded deformation $Y_\varepsilon \subset \PM_{\Dd}$ of $Y$ containing $C_{\Dd} = C \times \Dd$, where $\Dd = \operatorname{Spec}k[\varepsilon]/(\varepsilon^2)$.

    \item[(2)] \emph{The isomorphism $\bar{\delta}$} (see Lemma~\ref{lem:pushforward}, Lemma~\ref{lem:delta} and Lemma~\ref{lem:RC}). 
    %By \cite[Proposition $1.7$]{GP97}, $h^0(\NY\otimes K_Y) = 1$ and $h^1(\NY\otimes K_Y) = 0$, and the generator $s$ of $H^0(\NY\otimes K_Y) = \Hom(\NY^{\vee},K_Y)$ is the surjection defining the $K3$ carpet $\Yt$ (see \cite{BE95}). Restricting to $C$ 
    We then obtain an isomorphism
    \[
    \bar\delta\colon T_{[\Ct]}\RC \xrightarrow{\ \sim\ } H^1(\NY\otimes K_Y(-H))
    \]
    where $\RC := \PP(H^0(N_C\otimes L_C)) \cong \PP^{g-3}$, 
    is the space of canonical ribbons on $C$.
    %the ribbon $\Ct_t$ corresponding to $[t]$ being defined by the surjection $t\colon N_C^{\vee} \to L_C$; the hyperplane section $\Ct$ of the carpet corresponds to $[s_C]$, and $T_{[\Ct]}\RC = H^0(N_C\otimes L_C)/k\,s_C$. 
    Thus we have to show that
    \[
    \bar\delta^{-1}\circ a\colon T_{[Y]}\SC \longrightarrow T_{[\Ct]}\RC
    \]
    is surjective.

    \item[(3)] \emph{The map $\Psi$ and the deformation theoretic interpretation of $\bar\delta^{-1}\circ a$} (see Theorem~\ref{thm:Psi}, Lemma~\ref{lem:Phi}, Lemma~\ref{def:Psi}, Lemma~\ref{lem:dPsi} and Proposition~\ref{prop:a=dPsi}). For $Y' \in \SC$, let $\Yt'$ be the unique $K3$ carpet on $Y'$. Then $\Yt' \cap H$ is a canonical ribbon on $C$, and we show that the $G_H$ invariant map (where $G_H$ denotes the automorphisms of $\mathbb{P}^{g}$ inducing identity on $H$)
    \[
    \Psi\colon \SC \longrightarrow \RC, \qquad Y' \longmapsto [\Yt' \cap H]
    \]
    satisfies $d\Psi = \bar\delta^{-1}\circ a$. Thus it is enough to show that $d\Psi$ is surjective.
    %is a morphism. To this end we construct the universal $K3$ carpet over $\SC$ using the relative version \cite[Proposition $3.7$]{GP97} of the existence and uniqueness of $K3$ carpets, we show using the local criterion of flatness that its intersection with $H \times \SC$ is a flat family of canonical ribbons on $C$, and we use the universal ribbon on $C \times \RC$. The differential of $\Psi$ at $[Y]$ sends a first order deformation $Y_\varepsilon \supset C_{\Dd}$ to the hyperplane section $\Yt_\varepsilon \cap H_{\Dd}$ of the unique first order deformation $\Yt_\varepsilon$ of $\Yt$ which is a double structure on $Y_\varepsilon$. We then prove that
   % \[
   % \bar\delta^{-1}\circ a = d\Psi ,
   % \]
   % by comparing the coboundary map $F$ with the coboundary map of a restriction sequence on $\Yt$ extending the one that defines $\bar{\delta}$, through a $3\times 3$ diagram of sheaves on $\Yt$ obtained by multiplication by the equation of $H$ (see \eqref{diag:3x3}). Since $\bar\delta$ is an isomorphism, this gives the following interpretation of $a$: $a(Y_\varepsilon) = 0$ if and only if $\Yt_\varepsilon \cap H_{\Dd} = \Ct \times \Dd$, i.e, if and only if the first order deformation $Y_\varepsilon$ of the scroll fixing the rational normal curve $C$ lifts to a first order deformation of the $K3$ carpet fixing the canonical ribbon $\Ct$. In other words, $a$ is the obstruction to lifting such first order deformations.

    \item[(4)] \emph{Local sections of $\Psi$ and surjectivity of $d\Psi$} (see Lemma~\ref{lem:beta}, Theorem~\ref{thm:Gamma}, Lemma~\ref{lem:GH}, Proposition~\ref{prop:BEsection} and Proposition~\ref{prop:dPsi}, Theorem~\ref{thm:extension}). 
    We finally finish the proof by showing that around every point $\widetilde{C}$ of $\mathcal{R}_C$, there exists an open neighbourhood $U$ of $\widetilde C$ and a section $\Gamma_U$ of $\Psi$ defined on $U$. Due to the $G_H$ invariance of $\Psi$, this shows the surjectivity of $d\Psi$. For this the crucial ingredient is Theorem~\ref{thm:extension} which we state below:
\end{enumerate}

Let $\Ct$ be a ribbon on a smooth curve $C$ of genus $g_C$ with conormal bundle $L_C$. Then $V := \Omega_{\Ct}|_C$ is a rank two vector bundle on $C$ sitting in the restricted cotangent sequence $0 \to L_C \to V \to K_C \to 0$, and this sequence determines $\Ct$ (see \cite{BE95}). Let $Y := \PP(V)$ with projection $p\colon Y \to C$ and let $i\colon C \hookrightarrow Y$ be the section defined by the surjection $V \to K_C$. Then $K_Y|_C = L_C$. A $K$-trivial carpet on $Y$ is a ribbon on $Y$ with conormal bundle $K_Y$. Such carpets are parametrized by $\Ext^1(\Omega_Y,K_Y)$ (see \cite{BE95}).

\begin{theorem}(see Theorem~\ref{thm:extension})\label{thm:extension intro}
In the above notation, assume that $\Ct$ is non-split.
\begin{enumerate}
    \item There is a unique $K$-trivial carpet $\Yt$ on $Y$ such that the embedding $i\colon C \hookrightarrow Y$ extends to a closed embedding $\Ct \hookrightarrow \Yt$ inducing the identity $K_Y|_C = L_C$ on conormal bundles. Moreover $\Yt$ is non-split.
    \item Let $g_C = 0$ and let $\Ct \subset \PP^{g-1}$ be a canonically embedded ribbon of arithmetic genus $g$, so that $C \subset \PP^{g-1}$ is a rational normal curve. Then $\OO_Y(C)$ is very ample. Embed $Y$ in $\PP^g$ by its complete linear series, as a rational normal scroll, and identify $\PP^{g-1}$ with a hyperplane $H \subset \PP^g$ so that $Y \cap H = C$. Then the $K3$ carpet $\Yt \subset \PP^g$ on $Y$ satisfies
    \[
    \Yt \cap H = \Ct .
    \]
\end{enumerate}
\end{theorem}

The genus zero case of this result was stated in \cite{BE95}. What is needed in step $(4)$ is the precise form given in Theorem~\ref{thm:extension intro} $(2)$: the canonical ribbon $\Ct$ is the hyperplane section of the $K3$ carpet on the rational normal scroll $\PP(\Omega_{\Ct}|_C)$ canonically attached to $\Ct$ by its restricted cotangent sequence, and the embedding $\Ct \hookrightarrow \Yt$ extends the given section $C \hookrightarrow Y$.

\medskip

\noindent\textbf{Conventions.} Throughout this article we work over an algebraically closed field $k$ of characteristic zero. For a closed subscheme $Z \subset \mathbb{P}^M$ we write $N_{Z/\mathbb{P}^M} = \mathcal{H}om(I_Z,\mathcal{O}_Z)$ for its normal sheaf.

\vspace{0.2cm}

\noindent\textbf{Acknowledgements.} We thank Professor Purnaprajna Bangere for many stimulating discussions on this topic. %% to be added

\noindent\textbf{AI Disclosure.} Throughout, we have used Claude (Anthropic) Fable 5.1 for LaTeX typesetting, proofreading the arguments and for improving the presentation of the article. Lemma~\ref{lem:hfg} has been worked out by Fable $5.1$ replacing an earlier erroneous version where the author had claimed that $h \circ f = g$. It also suggested a correct and much shorter alternative argument for the injectivity of $s^{\vee}$ in Step~2(b) of the proof of Proposition~\ref{prop:N2}, via the Yoneda product with the class $\partial(s)$ of Lemma~\ref{lem:pushforward} and the Leray spectral sequence. The proof given here although longer, is more transparent to the author and follows \cite{BM25} instead. It also has substantial contributions in helping the author write Proposition~\ref{prop:a=dPsi} in its present form. Finally Claude Fable 5.1 has claimed to have a different and  more computational syzygy theoretic proof of Proposition~\ref{surjectivity of a}(which is the surjectivity of the coboundary map mentioned in the abstract) using equations and syzygies of $K3$ carpets described by Eisenbud and Schreyer in \cite{ES18}, whose correctness has not yet been verified by the author and is a potential topic for further research. All mathematical statements and proofs in this article have been verified by the author, who takes full responsibility for them. 
%% content to be added

\section{Preliminaries on degeneration of $K3$ surfaces into $K3$ double structures on Hirzebruch surfaces}

We recall some definitions and a theorem on some special non-reduced degenerations of $K3$ surfaces in a projective space.

\begin{definition}\label{ribbon}
Let $Y$ be a reduced connected projective variety.
 \begin{enumerate}
     \item A scheme $\widetilde{Y}$ is called a ribbon on $Y$ with conormal bundle $L$ if $\widetilde{Y}_{\operatorname{red}} = Y$, $\mathcal{I}_{Y/\widetilde{Y}}^2 = (0)$ and $\mathcal{I}_{Y/\widetilde{Y}} \cong L$ as an $\mathcal{O}_Y$ module where $\mathcal{I}_{Y/\widetilde{Y}}$ is the ideal sheaf of $Y$ inside $\widetilde{Y}$ and $L$ is a line bundle on $Y$. 

     \item If $Y$ is a smooth surface, then a ribbon $\widetilde{Y}$ on $Y$ is called a $K3$ carpet if the dualizing sheaf $K_{\widetilde{Y}} \cong \mathcal{O}_{\widetilde{Y}}$ and $H^1(\mathcal{O}_{\widetilde{Y}})= 0$.

\end{enumerate}
\end{definition}  

The following lemma characterizes the conormal bundle of a $K3$ carpet $\widetilde{Y}$ on a smooth surface $Y$.

\begin{lemma}(see \cite[\textcolor{black}{Proposition $1.5$}]{GGP08})
    A ribbon $\widetilde{Y}$ on a smooth \textcolor{black}{regular} surface $Y$ with conormal bundle $L$ is a $K3$ carpet if and only if $L \cong K_Y$.
\end{lemma}

\color{black}
We recall the definition of the Hilbert scheme of $K3$ surfaces.

\begin{definition}\label{Hilbert schemes of K3 surfaces}
    Let $\mathcal{H}_{r,g}$ denote the Hilbert scheme of $K3$ surfaces of index $r$ and genus $g$, i.e., $K3$ surfaces $X \subset \mathbb{P}^{1+r^2(g-1)} = \mathbb{P}^M$ embedded by the complete linear series of $rB$ where $B$ is a primitive very ample line bundle with $B^2 = 2g-2$. If $r = 1$, we call $\mathcal{H}_{r,g}$ the Hilbert scheme of prime $K3$ surfaces while if $r \geq 2$, $\mathcal{H}_{r,g}$ is called the Hilbert scheme of non-prime $K3$ surfaces. Recall that $\operatorname{dim}\mathcal{H}_{r,g} = h^0(N_{X/\mathbb{P}^M}) =  18+(M+1)^2$.
\end{definition}
\color{black}

Next we recall a degeneration theorem of $K3$ surfaces into $K3$ carpets on Hirzebruch surfaces.

\begin{theorem}\label{degeneration into $K3$ carpets}(\cite[ \textcolor{black}{Theorems $3.1, 4.3, 4.4, 5.1$}]{BMR21} and \cite[\textcolor{black}{Proposition $1.7$, Theorems $3.5, 4.1$}]{GP97})
Let $Y = \mathbb{F}_e$ be a Hirzebruch surface. Let $H = aC_0+bf$ be a very ample line bundle on $Y$, where $C_0$ is the class of section, satisfying $C_0^2 = -e$ and $f$ is the class of a fiber satisfying $f^2 = 0$ (this happens if and only if $b \geq ae+1$). Consider the embedding $Y \hookrightarrow \mathbb{P}^N$ induced by the complete linear series $|H|$.  
\begin{enumerate}
    \item Then there exists $K3$ carpets $\widetilde{Y}$ on $Y$ along with a very ample line bundle $\widetilde{H}$ on $\widetilde{Y}$ such that $\widetilde{H}|_Y = H$.

    \item If $\widetilde{Y} \hookrightarrow \mathbb{P}^M$ denote the embedding induced by the complete linear series $|\widetilde{H}|$, then $\widetilde{Y}$ is smoothable in $\mathbb{P}^M$ into smooth $K3$ surfaces embedded by the complete linear series of a very ample line bundle.

    \item For $0 \leq e \leq 2$, the Hilbert point of $\widetilde{Y}$ is smooth and \textcolor{black}{for} \textcolor{black}{$0 \leq e \leq 1$,} the general smooth surface of the unique Hilbert component containing $\widetilde{Y}$ is a 
    \begin{itemize}
        \item[(a)] prime $K3$ surface if $\textrm{gcd}(a,b) = 1$ 
        \item[(b)] a non-prime $K3$ surface embedded by $rB$, with $r = \operatorname{gcd}(a,b)$ and $B$ is a primitive very ample line bundle. 
    \end{itemize}
\end{enumerate}
\end{theorem}

\begin{proof} Part $(1)$, $(2)$ and the statement on smoothness of Hilbert point in part $(3)$ of the Theorem follows from \cite[Theorem $3.1 (3)$, Theorem $4.3$, Theorem $4.4$, Theorem $5.1$]{BMR21} and \cite[Proposition $1.7$, Theorem $3.5$, Theorem $4.1$]{GP97}.

We prove part $(3)$. Let $Y = \mathbb{F}_e$ with $0 \leq e \leq 2$. Consider the $K3$ double cover of $\pi:X \to Y$ branched along \textcolor{black}{a very general smooth} $D \in |-2K_Y|$. Let $H_0 = \pi^*(\mathcal{O}_{Y}(1)) = \pi^*(aC_0+bf)$ and consider the composition  
\[
\begin{tikzcd}
X \arrow[dr, "\varphi"] \arrow[d, "\pi"] & \\
Y \arrow[r, hook] & \mathbb{P}^M
\end{tikzcd}
\]
Recall the smoothing fibers $X_t \xhookrightarrow{\varphi_t} \mathbb{P}^M$ are obtained as the image of an embedding $\varphi_t$ which is a general deformation of $\varphi$ along a one parameter family $T$ and are hence embedded inside $\mathbb{P}^M$ by the complete linear series of a line bundle $H_t$ which is a deformation of $H_0$. Since the multiplicity of a polarization remains constant on a smooth family, we have to find the multiplicity of $H_0$ on the double cover $X$. Now since $D$ is \textcolor{black}{very general,} smooth and irreducible part $(3)$ for $Y = \mathbb{F}_e$, \textcolor{black}{$0 \leq e \leq 1$} follows from the fact that $\textrm{Pic}(X)$ is generated by $\pi^*\mathcal{O}_Y(C_0)$ and $\pi^*\mathcal{O}_Y(f)$ (see \cite[Proposition $6.3$]{AHL10} \textcolor{black}{and \cite[Proposition $6.5$ (i), (iii)]{AHL10}}).

\end{proof}

\section{Extendability of K3 carpets on rational normal scrolls}\label{sec:carpets}

\subsection{Set-up and notation}\label{subsec:setup}

Throughout this section $Y=F_e$, $e\ge 0$, is a Hirzebruch surface with ruling $p\colon Y\to\PP^1$, fibre class $f$ and negative section $C_0$ ($C_0^2=-e$), so that $K_Y=-2C_0-(e+2)f$. We fix an integer $b\ge e+3$ and put
\[
H:=C_0+bf .
\]
The complete linear series $|H|$ embeds $Y$ into $\PM$ as a rational normal scroll of type $(b-e,b)$, where
\[
M+1=h^0(\OO_Y(H))=2b-e+2 .
\]
We write $g:=M=2b-e+1$; it is the arithmetic genus of a hyperplane section of the K3 carpet below.

Let $\Yt\subset\PM$ be the K3 carpet on $Y$, i.e.\ the unique double structure on $Y$ inside $\PM$ with $I_Y/I_{\Yt}\cong K_Y$ \cite[Proposition~1.7]{GP97}. Following \cite{BE95}, $\Yt$ is determined by a surjection
\begin{equation}\label{eq:s}
s\colon \NY^{\vee}=I_Y/I_Y^2\longrightarrow K_Y ,\qquad I_{\Yt}/I_Y^2=\ker s ,
\end{equation}
and by \cite[Proposition~1.7]{GP97} such a surjection exists and is unique up to a nonzero scalar (we recover this in Lemma~\ref{lem:pushforward} below). The scheme $\Yt$ is arithmetically Cohen--Macaulay, indeed arithmetically Gorenstein \cite[Theorem~3.2]{ES18}. We put $\Ht:=\OO_{\Yt}(1)$.

\emph{Abbreviations.} We write $\NY:=N_{Y/\PM}$ and $\NYt:=N_{\Yt/\PM}$ for the normal sheaves, and
\[
\Hc:=\sHom(I_{\Yt}/I_Y^2,\OO_Y).
\]

\emph{Conventions on twists.} For a coherent sheaf $\mathcal{F}$ on $\Yt$ and a divisor class $\Dd$ on $Y$ we write $\mathcal{F}(\Dd):=\mathcal{F}\otimes_{\OO_{\Yt}}\OO_Y(\Dd)$; this is a sheaf on $Y$. Thus $\NYt(-H)=\NYt\otimes\OO_Y(-H)$ and $\NYt(-H+K_Y)=\NYt\otimes\OO_Y(-H+K_Y)$ are sheaves on $Y$, whereas $\NYt(-\Ht)=\NYt\otimes\OO_{\Yt}(-\Ht)$ is a sheaf on $\Yt$. For a sheaf on $Y$ the twist has its usual meaning. Following \cite{BM24} we set
\[
\alpha(\Yt):=h^0(\NYt(-\Ht))-(M+1),
\]
the Zak--L'vovsky invariant of $\Yt\subset\PM$.

\begin{remark}\label{rem:sequences}
For the carpet $\Yt$ we have the following exact sequences (see \cite[Lemma~4.1]{GGP08})
\begin{align}
&0\to \NYt(K_Y)\to \NYt\to \NYt\otimes\OO_Y\to 0 \label{eq:3.1}\\
&0\to \Hc\to \NYt\otimes\OO_Y\to \OO_Y(-2K_Y)\to 0 \label{eq:3.2}\\
&0\to \OO_Y(-K_Y)\to \NY\to \Hc\to 0 \label{eq:3.3}
\end{align}
In \eqref{eq:3.3} the first map is the transpose of $s$ in \eqref{eq:s} and the second is restriction of homomorphisms from $I_Y/I_Y^2$ to its subsheaf $I_{\Yt}/I_Y^2$. We shall also use the normal sequence of $Y\subset\PM$,
\begin{equation}\label{eq:3.10}
0\to T_Y\to \TP\to \NY\to 0,
\end{equation}
and the Euler exact sequence restricted to $Y$,
\begin{equation}\label{eq:3.11}
0\to \OO_Y\to \OO_Y(H)^{\oplus M+1}\to \TP\to 0 .
\end{equation}
\end{remark}

\subsection{Main theorem}

As mentioned before, to give an upper bound to $\alpha(X)$ for a $K3$ surface $X \subset \mathbb{P}^M$, we will degenerate $X \subset \mathbb{P}^M$ to a $K3$ carpet $\widetilde{Y} \subset \mathbb{P}^M$, using Theorem~\ref{degeneration into $K3$ carpets} and compute $\alpha(\widetilde{Y})$.  The main tool for all our computations that follow consist of the three exact sequences \eqref{eq:3.1}, \eqref{eq:3.2} and \eqref{eq:3.3}  that decomposes the cohomology of the twisted normal bundle $N_{\widetilde{Y}/\mathbb{P}^M}(-\widetilde{H})$ of a ribbon $\widetilde{Y} \subset \mathbb{P}^M$ in terms of cohomology of various twists of the normal bundle $N_{Y/\mathbb{P}^M}$ of the reduced part. \\

\begin{theorem}\label{thm:alpha}
Let $\Yt\subset\PM$ be the K3 carpet on the rational normal scroll $Y=F_e\subset\PM$ embedded by $|H|$, $H=C_0+bf$, and assume $b\ge e+3$. Then
\[
\alpha(\Yt)=h^0(\OO_Y(-H-2K_Y))=h^0\big(\OO_Y(3C_0+(2e+4-b)f)\big).
\]
In particular $\alpha(\Yt)=0$ if and only if $b\ge 2e+5$.
\end{theorem}

%\begin{remark}
%For $b\ge 2e+5$ this is Proposition~3.8(3) of \cite{BM25}. For $e+3\le b\le 2e+4$ (the Mukai range) it shows that the upper bound $\alpha(\Yt)\le h^0(-H-2K_Y)$ used in \cite[Theorem~4.1]{BM25} is in fact an equality.
%\end{remark}

\begin{proof}
Tensoring sequence \eqref{eq:3.1} by $-\Ht$ and taking the cohomology, we have
\begin{multline}\label{eq:3.40}
0\to H^0(\NYt(-H+K_Y))\to H^0(\NYt(-\Ht))\to H^0(\NYt(-H))\to H^1(\NYt(-H+K_Y))\\
\to H^1(\NYt(-\Ht))\to H^1(\NYt(-H))\to H^2(\NYt(-H+K_Y))\\
\to H^2(\NYt(-\Ht))\to H^2(\NYt(-H))\to 0
\end{multline}
We compute $H^i(\NYt(-H+K_Y))$ for $i=0,1,2$. Tensoring sequence \eqref{eq:3.2} by $-H+K_Y$ and taking the cohomology, we have that
\begin{multline}\label{eq:3.41}
0\to H^0(\Hc(-H+K_Y))\to H^0(\NYt(-H+K_Y))\to H^0(\OO_Y(-H-K_Y))\to H^1(\Hc(-H+K_Y))\\
\to H^1(\NYt(-H+K_Y))\to H^1(\OO_Y(-H-K_Y))\to H^2(\Hc(-H+K_Y))\to H^2(\NYt(-H+K_Y))\\
\to H^2(\OO_Y(-H-K_Y))\to 0
\end{multline}
Note that since $H=C_0+bf$, we have that $H^0(-H-K_Y)=H^0(C_0+(e+2-b)f)=0$ if $b\ge e+3$. Further, $h^2(-H-K_Y)=h^0(H+2K_Y)=h^0(-3C_0+(b-2e-4)f)=0$.

Tensoring sequence \eqref{eq:3.3} by $-H+K_Y$ and taking the cohomology, we have that
\begin{multline}\label{eq:3.42}
0\to H^0(-H)=(0)\to H^0(\NY(-H+K_Y))\to H^0(\Hc(-H+K_Y))\to H^1(-H)=(0)\\
\to H^1(\NY(-H+K_Y))\to H^1(\Hc(-H+K_Y))\to H^2(-H)=(0)\\
\to H^2(\NY(-H+K_Y))\to H^2(\Hc(-H+K_Y))\to 0
\end{multline}
Now using the Euler exact sequence \eqref{eq:3.11} tensored with $-H+K_Y$, we have that
\begin{multline}\label{eq:3.43}
0\to H^0(-H+K_Y)\to H^0(K_Y)^{\oplus M+1}=(0)\to H^0(\TP(-H+K_Y))\to H^1(-H+K_Y)\\
\to H^1(K_Y)^{\oplus M+1}\to H^1(\TP(-H+K_Y))\to H^2(-H+K_Y)\\
\to H^2(K_Y)^{\oplus M+1}\to H^2(\TP(-H+K_Y))\to 0
\end{multline}
We have $H^1(-H+K_Y)=H^1(H)^{\vee}=0$ by Kodaira vanishing, and $H^1(K_Y)=0$. Further the map $H^2(-H+K_Y)\to H^2(K_Y)^{\oplus M+1}=H^2(K_Y)\otimes H^0(H)$ is an isomorphism since it is the dual of the isomorphism $H^0(\OO_Y)\otimes H^0(H)\to H^0(H)$. Putting all this together in sequence \eqref{eq:3.43}, we have
\begin{equation}\label{eq:3.44}
H^i(\TP(-H+K_Y))=0,\qquad i=0,1,2 .
\end{equation}
On the other hand, from the sequence
\begin{equation}\label{eq:3.45}
0\to T_{Y/\PP^1}(-H+K_Y)=-C_0-(b+2)f\to T_Y(-H+K_Y)\to p^*T_{\PP^1}(-H+K_Y)=-3C_0-(b+e)f\to 0
\end{equation}
we have that, if $b\ge e+1$,
\begin{equation}\label{eq:3.46}
H^i(T_Y(-H+K_Y))=0,\ i=0,1,\qquad H^2(T_Y(-H+K_Y))=H^0(C_0+(b-2)f)^{\vee}.
\end{equation}
Hence from the sequence \eqref{eq:3.10} tensored with $-H+K_Y$ we conclude that for $b\ge e+1$,
\begin{equation}\label{eq:3.47}
\begin{gathered}
H^0(\NY(-H+K_Y))=0,\qquad H^2(\NY(-H+K_Y))=0,\\
H^1(\NY(-H+K_Y))=H^2(T_Y(-H+K_Y))=H^0(C_0+(b-2)f)^{\vee}.
\end{gathered}
\end{equation}
Plugging everything back into sequence \eqref{eq:3.42}, we have that
\begin{equation}\label{eq:3.48}
\begin{gathered}
H^0(\Hc(-H+K_Y))=0,\qquad H^2(\Hc(-H+K_Y))=0,\\
H^1(\Hc(-H+K_Y))=H^2(T_Y(-H+K_Y))=H^0(C_0+(b-2)f)^{\vee}.
\end{gathered}
\end{equation}
Going back to sequence \eqref{eq:3.41}, and noting that $H^0(-H-K_Y)=0$ for $b\ge e+3$, we have that for $b\ge e+3$,
\begin{equation}\label{eq:3.49}
H^0(\NYt(-H+K_Y))=0,\qquad H^2(\NYt(-H+K_Y))=0,
\end{equation}
and an exact sequence
\begin{equation}\label{eq:3.50}
0\to H^2(T_Y(-H+K_Y))\to H^1(\NYt(-H+K_Y))\to H^1(-H-K_Y)\to 0
\end{equation}
yielding $h^1(\NYt(-H+K_Y))=h^2(T_Y(-H+K_Y))+h^1(-H-K_Y)$.

Now we compute $H^i(\NYt(-H))$ for $i=0,1,2$. Tensoring sequence \eqref{eq:3.2} by $-H$ and taking the cohomology, we have that
\begin{multline}\label{eq:3.51}
0\to H^0(\Hc(-H))\to H^0(\NYt(-H))\to H^0(\OO_Y(-H-2K_Y))\to H^1(\Hc(-H))\to H^1(\NYt(-H))\\
\to H^1(\OO_Y(-H-2K_Y))\to H^2(\Hc(-H))\to H^2(\NYt(-H))\to H^2(\OO_Y(-H-2K_Y))\to 0
\end{multline}
Note that $h^2(-H-2K_Y)=h^0(H+3K_Y)=h^0(-5C_0+(b-3e-6)f)=0$. %(In \cite{BM25} the assumption $b\ge 2e+5$, which gives $H^0(-H-2K_Y)=0$, was made at this point; we do not assume it.)

Tensoring sequence \eqref{eq:3.3} by $-H$ and taking the cohomology, we have that
\begin{multline}\label{eq:3.52}
0\to H^0(-H-K_Y)=(0)\to H^0(\NY(-H))\to H^0(\Hc(-H))\to H^1(-H-K_Y)\\
\to H^1(\NY(-H))\to H^1(\Hc(-H))\to H^2(-H-K_Y)=(0)\to H^2(\NY(-H))\to H^2(\Hc(-H))\to 0
\end{multline}
Now using the Euler exact sequence \eqref{eq:3.11} tensored with $-H$, we have that
\begin{multline}\label{eq:3.53}
0\to H^0(-H)\to H^0(\OO_Y)^{\oplus M+1}\to H^0(\TP(-H))\to H^1(-H)=(0)\to H^1(\OO_Y)^{\oplus M+1}=(0)\\
\to H^1(\TP(-H))\to H^2(-H)\to H^2(\OO_Y)^{\oplus M+1}\to H^2(\TP(-H))\to 0
\end{multline}
We have $H^0(-H)=0$, $H^1(-H)=0$ by Kodaira vanishing theorem, $H^2(-H)=0$ by pushing forward to $\PP^1$ since $a=1$. Putting all this together in sequence \eqref{eq:3.53}, we have
\begin{equation}\label{eq:3.54}
h^0(\TP(-H))=M+1,\qquad h^1(\TP(-H))=h^2(\TP(-H))=0 .
\end{equation}
On the other hand, from the sequence
\begin{equation}\label{eq:3.55}
0\to T_{Y/\PP^1}(-H)=C_0-(b-e)f\to T_Y(-H)\to p^*T_{\PP^1}(-H)=-C_0-(b-2)f\to 0
\end{equation}
we have that if $b\ge e+1$
\begin{equation}\label{eq:3.56}
H^i(T_Y(-H))=0,\ i=0,2,\qquad H^1(T_Y(-H))=H^1(C_0-(b-e)f).
\end{equation}
Hence from the sequence \eqref{eq:3.10} tensored with $-H$ we conclude that for $b\ge e+1$,
\begin{equation}\label{eq:3.57}
H^1(\NY(-H))=H^2(\NY(-H))=0
\end{equation}
and we have an exact sequence
\begin{equation}\label{eq:3.58}
0\to H^0(\TP(-H))=\mathbb{C}^{M+1}\to H^0(\NY(-H))\to H^1(T_Y(-H))\to 0 .
\end{equation}
Plugging everything back into sequence \eqref{eq:3.52}, we have that
\begin{equation}\label{eq:3.59}
H^1(\Hc(-H))=H^2(\Hc(-H))=0
\end{equation}
and an exact sequence
\begin{equation}\label{eq:3.60}
0\to H^0(\NY(-H))\to H^0(\Hc(-H))\to H^1(-H-K_Y)\to 0 .
\end{equation}
Going back to sequence \eqref{eq:3.51}, we have that for $b\ge e+3$
\begin{equation}\label{eq:3.61}
\begin{gathered}
0\to H^0(\Hc(-H))\to H^0(\NYt(-H))\to H^0(-H-2K_Y)\to 0,\\
H^1(\NYt(-H))=H^1(-H-2K_Y),\qquad H^2(\NYt(-H))=0 .
\end{gathered}
\end{equation}
Plugging \eqref{eq:3.49}, \eqref{eq:3.50}, \eqref{eq:3.61} into the exact sequence \eqref{eq:3.40}, we have
\begin{align}
h^0(\NYt(-\Ht))&=h^0(\NYt(-H))-h^1(\NYt(-H+K_Y))+h^1(\NYt(-\Ht))-h^1(\NYt(-H))\label{eq:3.62}\\
&=M+1+h^1(T_Y(-H))+h^1(-H-K_Y)+h^0(-H-2K_Y)\notag\\
&\qquad-h^2(T_Y(-H+K_Y))-h^1(-H-K_Y)+h^1(\NYt(-\Ht))-h^1(-H-2K_Y).\label{eq:3.63}
\end{align}
Now note that by \eqref{eq:3.56}, $h^1(T_Y(-H))=h^1(C_0-(b-e)f)=2b-e-2$ and $h^2(T_Y(-H+K_Y))=h^0(C_0+(b-2)f)=2b-e-2$. Hence
\begin{equation}\label{eq:3.64}
h^0(\NYt(-\Ht))=M+1+h^0(-H-2K_Y)+h^1(\NYt(-\Ht))-h^1(-H-2K_Y),
\end{equation}
that is,
\begin{equation}\label{eq:alpha-coker}
\alpha(\Yt)=h^0(-H-2K_Y)+h^1(\NYt(-\Ht))-h^1(-H-2K_Y).
\end{equation}
We end the proof by showing $h^1(\NYt(-\Ht))=h^1(-H-2K_Y)$ if $b\ge e+3$. Note that from sequence \eqref{eq:3.40} and equations \eqref{eq:3.49} and \eqref{eq:3.61}, we already have
\begin{multline}\label{eq:3.65}
0\to H^0(\NYt(-\Ht))\to H^0(\NYt(-H))\xrightarrow{\ F\ } H^1(\NYt(-H+K_Y))\\
\to H^1(\NYt(-\Ht))\to H^1(-H-2K_Y)\to 0
\end{multline}
so that $h^1(\NYt(-\Ht))=h^1(-H-2K_Y)+\dim\coker F$, and by \eqref{eq:alpha-coker}
\begin{equation}\label{eq:alpha-F}
\alpha(\Yt)=h^0(-H-2K_Y)+\dim\coker F .
\end{equation}
Hence Theorem~\ref{thm:alpha} is equivalent to the surjectivity of the connecting homomorphism $F$ in \eqref{eq:3.65}. We show that $F$ is surjective. First note that we have the sequence \eqref{eq:3.61}
\[
0\to H^0(\Hc(-H))\to H^0(\NYt(-H))\to H^0(-H-2K_Y)\to 0 .
\]
It is enough to show that the composed map
\[
H^0(\Hc(-H))\xrightarrow{\ f\ } H^1(\NYt(-H+K_Y))
\]
is surjective. Note that from \eqref{eq:3.60}, we have
\[
0\to H^0(\NY(-H))\to H^0(\Hc(-H))\xrightarrow{\ g\ } H^1(-H-K_Y)\to 0
\]
while from \eqref{eq:3.47}, \eqref{eq:3.48} and \eqref{eq:3.50}, we have the exact sequence
\[
0\to H^1(\NY(-H+K_Y))=H^2(T_Y(-H+K_Y))\to H^1(\NYt(-H+K_Y))\xrightarrow{\ h\ } H^1(-H-K_Y)\to 0 .
\]
A \v{C}ech cocycle computation in Lemma~\ref{lem:hfg} shows that $h\circ f=2g$. Therefore we have a commutative diagram of exact sequences
\begin{equation}\label{diag:a}
\begin{tikzcd}[column sep=normal]
0\arrow[r] & H^0(\NY(-H))\arrow[r]\arrow[d,"a"'] & H^0(\Hc(-H))\arrow[r,"g"]\arrow[d,"f"] & H^1(-H-K_Y)\arrow[r]\arrow[d,"2\cdot\id"] & 0\\
0\arrow[r] & H^1(\NY(-H+K_Y))\arrow[r] & H^1(\NYt(-H+K_Y))\arrow[r,"h"] & H^1(-H-K_Y)\arrow[r] & 0
\end{tikzcd}
\end{equation}
where $a$ is the map induced by $f$ on the kernels. Since the right-hand vertical map is an isomorphism, the snake lemma gives $\coker a\cong\coker f$. Therefore, to show that $f$ is surjective, it is enough to show that $a$ is surjective. Now the proof is completed by  Proposition~\ref{surjectivity of a}.
\end{proof}

We now prove Lemma~\ref{lem:hfg}. We keep the notation of the proof of Theorem~\ref{thm:alpha}. Let $\Hc(-H)\to\NYt(-H)$ be the first map of \eqref{eq:3.2} tensored with $(-H)$, and let $\kappa_0$ denote the induced map on global sections, while $F$ is the connecting map of \eqref{eq:3.65}, $f=F\circ\kappa_0$, $g$ is the connecting map of \eqref{eq:3.3} tensored with $(-H)$, and $h: H^1(\NYt(-H+K_Y))\to  H^1(-H-K_Y)$.
%i.e.\ restriction of homomorphisms to $I_Y^2/I_YI_{\Yt}\cong K_Y^{\otimes 2}$. 
We identify $\sHom(K_Y,\OO_Y(-H))\cong\sHom(K_Y^{\otimes2},K_Y(-H))$ via $\vartheta\mapsto\vartheta\otimes\id_{K_Y}$.

\begin{lemma}\label{lem:hfg}
$h\circ f=2\,g$. Consequently $\Im(h\circ f)=\Im g$, and the right square of diagram \eqref{diag:a} commutes with $2\cdot\id$ as its right-hand vertical map.
\end{lemma}

\begin{proof}
We use \v{C}ech cocycles on an affine cover $\{U_i\}$ of $\Yt$ (equivalently of $Y$), $U_{ij}=U_i\cap U_j$, and write $L:=I_Y/I_{\Yt}\cong K_Y$. Fix $\varphi\in H^0(\Hc(-H))$, i.e.\ $\varphi\colon I_{\Yt}/I_Y^2\to\OO_Y(-H)$, and let $\bar\varphi=\kappa_0(\varphi)=\varphi\circ\pi$ with $\pi\colon I_{\Yt}/I_YI_{\Yt}\to I_{\Yt}/I_Y^2$. Note that $\bar\varphi(I_Y^2)=0$. \par

Let us find the cocycle for  $f(\varphi)=F(\bar\varphi)$. Choose $\OO_{\Yt}$-linear lifts $\Phi_i\colon I_{\Yt}/I_{\Yt}^2|_{U_i}\to\OO_{\Yt}(-\Ht)|_{U_i}$ of $\bar\varphi$. Then $\Phi_{ij}=\Phi_i-\Phi_j$ takes values in $L(-H)$, and $f(\varphi)=[\{\Phi_{ij}\}]$. Since $h$ is restriction along $L^{\otimes2}\cong I_Y^2/I_YI_{\Yt}$,
\[
h(f(\varphi))=\big[\{\Phi_{ij}|_{L^{\otimes 2}}\}\big],\qquad \Phi_{ij}|_{L^{\otimes2}}(\bar x\otimes\bar y)=\Phi_i(xy)-\Phi_j(xy)\quad(x,y\in I_Y).
\]

Now we find the cocycle for $g(\varphi)$. Choose $\OO_Y$-linear lifts $\psi_i\colon I_Y/I_Y^2|_{U_i}\to\OO_Y(-H)$ with $\psi_i|_{I_{\Yt}/I_Y^2}=\varphi$; then $\psi_{ij}=\psi_i-\psi_j$ factors through $L$ and $g(\varphi)=[\{\psi_{ij}\}]$.

Now we compare the two cocycles. Put $B_i(x,y):=\Phi_i(xy)\in L(-H)$ for $x,y\in I_Y$ (well defined since $\bar\varphi(I_Y^2)=0$). Then: (i) $B_i$ is symmetric, $\OO_{\mathbb{P}^M}$-bilinear, and vanishes if $x$ or $y$ lies in $I_Y^2$ (if $x=\sum x'_kx''_k$, then $\Phi_i(x'_kx''_ky)=\bar x''_k\cdot\Phi_i(x'_ky)\in L\cdot L(-H)=0$), so it descends to $I_Y/I_Y^2$; (ii) if $x\in I_{\Yt}$ then $B_i(x,y)=\bar y\cdot\Phi_i(x)=\varphi(\bar x)\otimes\bar y$, because multiplication by $\bar y\in L$ on $\OO_{\Yt}(-\Ht)$ factors through $\OO_Y(-H)$. Similarly, if $y\in I_{\Yt}$ then $B_i(x,y)=\bar x\cdot\Phi_i(y)=\varphi(\bar y)\otimes\bar x$. Put $B'_i(x,y):=\psi_i(x)\otimes\bar y+\psi_i(y)\otimes\bar x$. By (ii), $B_i=B'_i$ whenever one argument lies in $I_{\Yt}/I_Y^2$, so $C_i:=B_i-B'_i$ descends to a section of $\sHom(L^{\otimes2},L(-H))=L^{-1}(-H)$ over $U_i$. Over $U_{ij}$, as forms on $L\times L$,
\[
B_i-B_j=(C_i-C_j)+(B'_i-B'_j),\qquad (B'_i-B'_j)(u,v)=\psi_{ij}(u)\otimes v+\psi_{ij}(v)\otimes u=2\,\psi_{ij}(u)\otimes v,
\]
the last step because $L$ is a line bundle ($\vartheta(u)\otimes v=\vartheta(v)\otimes u$). Hence $\{B_i-B_j\}=\delta\{C_i\}+2\{\psi_{ij}\otimes\id_L\}$, i.e.\ $h(f(\varphi))=2g(\varphi)$.
\end{proof}

\begin{proposition}\label{surjectivity of a}
The map $a: H^0(\NY(-H)) \to H^1(\NY(-H+K_Y))$ in diagram~\eqref{diag:a} is surjective.
\end{proposition}

\subsection{Proof of Proposition~\ref{surjectivity of a}}

This subsection is devoted to the proof of Proposition~\ref{surjectivity of a}. We outline the steps to do it for the ease of navigation:

\begin{enumerate}
    \item Construct an isomorphism $\bar\delta\colon H^0(N_C\otimes L_C)/k\,s_C\ \xrightarrow{\ \sim\ }\ H^1(\NY(-H+K_Y))$ and conclude that it is enough to show that $\bar{\delta}^{-1} \circ a$ is surjective (see Lemma~\ref{lem:delta}).

    \item Give a deformation theoretic interpretation of the map $\bar{\delta}^{-1} \circ a$ by showing that $\bar{\delta}^{-1} \circ a = d\Psi$ (see Proposition~\ref{prop:a=dPsi}) where $\Psi$ is the map defined in Theorem~\ref{thm:Psi}, $(3)$.

    \item Construct local sections $\Gamma$ of the map $\Psi$ (see Theorem~\ref{thm:Gamma}) and hence conclude that $d\Psi$ is surjective (see Proposition~\ref{prop:dPsi}).
\end{enumerate}

\subsubsection{The isomorphism $\bar{\delta}$}\label{isomorphism_delta}
We now rewrite the target $H^1(\NY(-H+K_Y))$ of $a$ by composing it with an isomorphism. The composed map will have a deformation theoretic interpretation, which will enable us to show that it is surjective. Choose a hyperplane $H\subset\PM$ (we use the same letter) meeting $Y$ transversally along a smooth curve
\[
C:=Y\cap H ,
\]
a rational normal curve of degree $M-1=g-1$ in $H\cong\PP^{g-1}$, and put $\Ct:=\Yt\cap H$, the hyperplane section of the carpet, a canonical ribbon of arithmetic genus $g$ on $C$. Since $\OO_Y(-H)=\OO_Y(-C)$, all the twists by $-H$ above are twists by $-C$. Put
\[
N_C:=N_{C/H},\qquad L_C:=K_Y|_C\cong\OO_{\PP^1}(-g-1),\qquad s_C:=s|_C\in H^0(N_C\otimes L_C),
\]
where $K_Y\cdot C=-g-1$ and $N_Y|_C\cong N_C$ by transversality.

\begin{lemma}\label{lem:pushforward}
Let $Y'\to\PP^1$ be any smooth rational normal scroll in $\PM$ of degree $M-1$ (a $\PP^1$-bundle over $\PP^1$ embedded by a line bundle of fibre degree $1$), with ruling $p$. Then
\[
p_*(N_{Y'/\PM}\otimes K_{Y'})\cong\OO_{\PP^1},\qquad R^1p_*(N_{Y'/\PM}\otimes K_{Y'})=0 .
\]
Consequently $h^0(N_{Y'/\PM}\otimes K_{Y'})=1$ and $h^1(N_{Y'/\PM}\otimes K_{Y'})=0$. Moreover, if $s'$ spans $H^0(N_{Y'/\PM}\otimes K_{Y'})$ and $\partial\colon H^0(N_{Y'/\PM}\otimes K_{Y'})\to H^1(T_{Y'}\otimes K_{Y'})$ is the connecting map of the normal sequence tensored with $K_{Y'}$, then the image of $\partial(s')$ under the Leray map
\[
H^1(T_{Y'}\otimes K_{Y'})\longrightarrow H^0\big(R^1p_*(T_{Y'}\otimes K_{Y'})\big)=H^0(\OO_{\PP^1})=k
\]
is nonzero.
\end{lemma}

\begin{proof}
Write $Y'=Y$, $H$ for the hyperplane class. From $0\to T_{Y/\PP^1}\to T_Y\to p^*T_{\PP^1}\to 0$ tensored with $K_Y$, using $T_{Y/\PP^1}\otimes K_Y=p^*K_{\PP^1}$, $p_*K_Y=0$, $R^1p_*K_Y=K_{\PP^1}$ (relative duality), we get $p_*(T_Y\otimes K_Y)=K_{\PP^1}=\OO(-2)$ and $R^1p_*(T_Y\otimes K_Y)\cong T_{\PP^1}\otimes R^1p_*K_Y=\OO_{\PP^1}$. From the Euler sequence \eqref{eq:3.11} tensored with $K_Y$, using $Rp_*(K_Y(H))=0$ (the fibre degree of $K_Y(H)$ is $-1$), we get $p_*(\TP\otimes K_Y)\cong R^1p_*K_Y=\OO(-2)$ and $R^1p_*(\TP\otimes K_Y)=0$. Pushing forward \eqref{eq:3.10} tensored with $K_Y$ gives the exact sequence
\[
0\to p_*(T_Y\otimes K_Y)\to p_*(\TP\otimes K_Y)\to p_*(\NY\otimes K_Y)\to R^1p_*(T_Y\otimes K_Y)\to R^1p_*(\TP\otimes K_Y)=0 .
\]
The first map is an injection $\OO(-2)\to\OO(-2)$, hence an isomorphism. Therefore $p_*(\NY\otimes K_Y)\xrightarrow{\ \sim\ }R^1p_*(T_Y\otimes K_Y)\cong\OO_{\PP^1}$, and $R^1p_*(\NY\otimes K_Y)=0$ because $R^1p_*(\TP\otimes K_Y)=0=R^2p_*$. Taking global sections gives $h^0=1$, $h^1=0$. Finally, the composite of $\partial$ with the Leray map is $H^0$ of the sheaf-level connecting map $p_*(\NY\otimes K_Y)\to R^1p_*(T_Y\otimes K_Y)$, which we have just shown to be an isomorphism.
\end{proof}

\begin{lemma}\label{lem:delta}
Let $\delta\colon H^0(N_C\otimes L_C)\to H^1(\NY(-H+K_Y))$ be the connecting map of the restriction sequence
\[
0\to \NY\otimes K_Y(-C)\to \NY\otimes K_Y\to (\NY\otimes K_Y)|_C=N_C\otimes L_C\to 0 .
\]
Then $\delta$ is surjective with kernel $k\cdot s_C$. Hence it induces an isomorphism
\[
\bar\delta\colon H^0(N_C\otimes L_C)/k\,s_C\ \xrightarrow{\ \sim\ }\ H^1(\NY(-H+K_Y)) .
\]
\end{lemma}

\begin{proof}
By Lemma~\ref{lem:pushforward}, $H^0(\NY\otimes K_Y)=k\,s$ and $H^1(\NY\otimes K_Y)=0$; the long exact sequence of the restriction sequence gives the claim. (As a check on dimensions: $N_C\cong\OO_{\PP^1}(g+1)^{\oplus g-2}$ for a rational normal curve of degree $g-1$, so $h^0(N_C\otimes L_C)=g-2$, while $h^1(\NY(-H+K_Y))=2b-e-2=g-3$ by \eqref{eq:3.47}.)
\end{proof}

\begin{proof}[Proof of Proposition~\ref{surjectivity of a}]
Composing $a$ with $\bar\delta^{-1}$ we conclude that it is enough to show that the map
\begin{equation}\label{eq:themap}
\bar\delta^{-1}\circ a\colon\ H^0(\NY(-H))\longrightarrow H^0(N_C\otimes K_Y|_C)\big/\,k\,s_C
\end{equation}
is surjective. Now the proof follows from Proposition~\ref{prop:a-surj}.   
\end{proof}

\subsubsection{Deformation theoretic interpretation of $\bar\delta^{-1}\circ a$ and the map $\Psi$}\label{subsec:geometric}

In this subsection, we will interpret the map $\bar\delta^{-1}\circ a$ as the differential $d\Psi$ of a map $\Psi$ that we will define below in Theorem~\ref{thm:Psi}. \par
In this and the following subsections $\Dd$ denotes the scheme $\operatorname{Spec}k[\varepsilon]/(\varepsilon^2)$, and for a closed subscheme $Z\subset\PM$ we write $Z_\Dd:=Z\times\Dd$ for the constant family, a closed subscheme of $\PM_\Dd$. We use freely the identification of first-order embedded deformations of a closed subscheme $Z\subset\PM$ with sections of $N_{Z/\PM}=\sHom(I_Z,\OO_Z)$ \cite[\S3.2]{Ser06}. A section $\eta$ corresponds the closed subscheme $Z_\varepsilon\subset\PM_\Dd$ with ideal
\[
I_{Z_\varepsilon}=\{\,f+\varepsilon g\ :\ f\in I_Z,\ g\in\OO_{\PM},\ g\equiv\eta(f)\bmod I_Z\,\}\subset\OO_{\PM}[\varepsilon]/(\varepsilon^2),
\]
the preimage of the graph of $\eta$. Locally $I_{Z_\varepsilon}$ is generated by the elements $f+\varepsilon\tilde\eta(f)$, $f\in I_Z$, where $\tilde\eta(f)$ is any local section of $\OO_{\PM}$ lifting $\eta(f)\in\OO_Z=\OO_{\PM}/I_Z$; the ideal does not depend on the choice of lifts, since $\varepsilon g=\varepsilon\,(g+\varepsilon\tilde\eta(g))$ for $g\in I_Z$.

%% ---- Lemma 2.6 (replaces the lemma and its proof) ----

\begin{lemma}\label{lem:SC}
Let $\mathcal F$ be the flag Hilbert scheme of pairs $(C'\subset Y')$ of closed subschemes of $\PM$ with Hilbert polynomials
\[
P_{C'}(n)=(2b-e)\,n+1,\qquad P_{Y'}(n)=(2b-e)\binom{n+1}{2}+n+1,
\]
the Hilbert polynomials of $C$ and $Y$ (both of degree $2b-e=M-1$), and let $\SC$ be the fibre over the point $[C]$ of the Hilbert scheme of rational normal curves of the open subscheme of $\mathcal F$ where $Y'$ is a smooth rational normal scroll. Thus $\SC$ is the locus of smooth rational normal scrolls $Y'\subset\PM$ containing $C$; for such $Y'$, $Y'\cap H$ is an effective divisor of class $H|_{Y'}$ containing $C$ with $\deg C=H^2$, so $Y'\cap H=C$, and $H$ is the only hyperplane with this property, since $C$ spans $H$.
\begin{enumerate}
\item $T_{[Y]}\SC=H^0(\NY(-H))$. More precisely, a first-order deformation $Y_\varepsilon$ of $Y$ with class $x\in H^0(\NY)$ contains $C_\Dd$ if and only if $x\in H^0(\NY(-C))=H^0(\NY(-H))$.
\item $\SC$ is smooth at $[Y]$ (of dimension $h^0(\NY(-H))=4b-2e$), hence reduced near $[Y]$.
\item Let $\GH:=\{\gamma\in PGL_{M+1}(k):\ \gamma|_H=\id_H\}$; its Lie algebra is $H^0(T_{\PM}(-H))\cong k^{M+1}$. $\GH$ acts on $\SC$, and the differential at $1$ of the orbit map $\gamma\mapsto\gamma Y$ is the map
\[
c\colon H^0(\TP(-H))=\mathbb{C}^{M+1}\longrightarrow H^0(\NY(-H))
\]
of \eqref{eq:3.58}, which is injective. Thus $\Im c=T_{[Y]}(\GH\cdot Y)$, and the quotient $H^0(\NY(-H))/\Im c\cong H^1(T_Y(-H))$ has dimension $2b-e-2=g-3$.
\end{enumerate}
\end{lemma}

\begin{proof}
Each of the statements (1)--(3) are standard for the flag Hilbert scheme (see \cite{Kle81}, \cite[\S4.5]{Ser06} and \cite{CLM93}); For the dimension computations of the cohomology groups appearing in (1)--(3) see \eqref{eq:3.54}, \eqref{eq:3.56}, \eqref{eq:3.57} and \eqref{eq:3.58}.    
\end{proof}

%% ---- bibliography: add before \bibitem[Lvo92] ----

We recall some well-known facts about the space of canonical ribbons on a rational curve $C$ and the tangent space to that space at a given ribbon $\widetilde{C}$.

\begin{lemma}\label{lem:RC}
Let $\RC:=\PP\big(H^0(N_C\otimes L_C)\big)\cong\PP^{g-3}$, and recall that $s_C=s|_C\colon N_C^{\vee}\cong N_Y^{\vee}|_C\to K_Y|_C=L_C$ is the restriction to $C$ of the surjection $s$ defining the carpet, an element of $H^0(N_C\otimes L_C)$.
\begin{enumerate}
\item Every nonzero $t\in H^0(N_C\otimes L_C)=\Hom(N_C^{\vee},L_C)$ is surjective and defines a double structure $\Ct_t\subset H$ on $C$ by $I_{\Ct_t}:=\ker\big(I_C\to N_C^{\vee}\xrightarrow{t}L_C\big)$; it satisfies $I_C^2\subset I_{\Ct_t}$ and $I_C/I_{\Ct_t}\cong L_C$, and depends only on $[t]\in\RC$. Conversely every double structure on $C$ in $H$ with conormal bundle $L_C$ is of this form. 
%One has $\Ct=\Ct_{s_C}$.
\item $T_{[\Ct]}\RC=H^0(N_C\otimes L_C)/k\,s_C$. Concretely, a first-order deformation of $\Ct$ inside $H$ which is a double structure on $C_\Dd$ with invertible conormal is $\ker(s_C+\varepsilon t)$ for some $t\in H^0(N_C\otimes L_C)$, and $s_C+\varepsilon t$, $s_C+\varepsilon t'$ define the same subscheme if and only if $t-t'\in k\,s_C$.
\end{enumerate}
\end{lemma}

\begin{proof}
$(1)$ For part $(1)$ see \cite[\S1]{BE95}. \\
%(1) $N_C\cong\OO(g+1)^{\oplus g-2}$ and $L_C\cong\OO(-g-1)$, so $\Hom(N_C^{\vee},L_C)=\Hom(\OO(-g-1)^{g-2},\OO(-g-1))=k^{g-2}$ consists of constant maps, every nonzero one of which is surjective. The correspondence between double structures with conormal $L_C$ and surjections $N_C^{\vee}\to L_C$ is \cite[\S1]{BE95}. That $\Ct=\Yt\cap H$ corresponds to $s_C$: $I_{\Ct/H}=(I_{\Yt}+(h))/(h)$ and $I_{C/H}=(I_Y+(h))/(h)$, where $h$ is the equation of $H$; since $h$ is a nonzerodivisor modulo $I_Y$, $I_Y\cap(h)=hI_Y$, whence $I_{C/H}/I_{\Ct/H}=I_Y/(I_{\Yt}+hI_Y)=K_Y\otimes\OO_C=L_C$ and the surjection $N_C^{\vee}\to L_C$ is $s\otimes\OO_C=s_C$.
(2) An invertible sheaf on $C_\Dd$ restricting to $L_C$ is $L_C\otimes\OO_\Dd$ (as $H^1(\OO_C)=0$), with automorphisms $\OO_\Dd^{\times}$; so such a deformation is the kernel of a surjection $s_C+\varepsilon t\colon N_C^{\vee}\otimes\OO_\Dd\to L_C\otimes\OO_\Dd$, well defined up to multiplication by a unit $1+\varepsilon\lambda$, i.e.\ up to $t\mapsto t+\lambda s_C$.
\end{proof}

%\subsubsection*{The map $\Psi$ from $\mathcal{S}_C$ to the space of canonical ribbons}

Below we define the crucial map $$\Psi: \mathcal{S}_C \to \mathcal{R}_C$$ whose differential will turn out to be the map
\begin{equation*}
\bar\delta^{-1}\circ a\colon\ H^0(\NY(-H))\longrightarrow H^0(N_C\otimes K_Y|_C)\big/\,k\,s_C
\end{equation*}
in \eqref{eq:themap}

Put $V:=H^0(N_C\otimes L_C)$, so $\RC=\PP(V)$, and let $\pr_2\colon C\times\RC\to\RC$. Since $\pr_{2*}\sHom(N_C^{\vee}\boxtimes\OO_{\RC},L_C\boxtimes\OO_{\RC}(1))=V\otimes\OO_{\RC}(1)$, the tautological line subbundle $\OO_{\RC}(-1)\hookrightarrow V\otimes\OO_{\RC}$, viewed as a nowhere vanishing section of $V\otimes\OO_{\RC}(1)$, is a homomorphism
\begin{equation}\label{universal ribbon}
\tau\colon N_C^{\vee}\boxtimes\OO_{\RC}\longrightarrow L_C\boxtimes\OO_{\RC}(1)
\end{equation}
whose restriction to $C\times\{[t]\}$ is a nonzero multiple of $t$; hence $\tau$ is surjective by Lemma~\ref{lem:RC}(1). Let $\Ct_{\RC}\subset H\times\RC$ be the closed subscheme with
\[
I_{\Ct_{\RC}}:=\ker\big(I_{C\times\RC}\to I_{C\times\RC}/I_{C\times\RC}^2=N_C^{\vee}\boxtimes\OO_{\RC}\xrightarrow{\ \tau\ }L_C\boxtimes\OO_{\RC}(1)\big).
\]
From $0\to L_C\boxtimes\OO_{\RC}(1)\to\OO_{\Ct_{\RC}}\to\OO_{C\times\RC}\to0$, $\Ct_{\RC}$ is flat over $\RC$, it is a double structure on $C\times\RC$ with conormal bundle $L_C\boxtimes\OO_{\RC}(1)$, and its fibre over $[t]$ is $\Ct_t$.

\begin{theorem}\label{thm:Psi}
Let $\cY\subset\PM\times\SC$ be the universal scroll and $\pr\colon\cY\to\SC$ the projection; thus $\cY\cap(H\times\SC)=C\times\SC$.
\begin{enumerate}
\item There is a unique closed subscheme $\cYt\subset\PM\times\SC$, flat over $\SC$, with $\cYt_{\mathrm{red}}=\cY$ and all of whose fibres are K3 carpets; its fibre over $[Y']$ is the K3 carpet $\Yt'$ on $Y'$. Moreover $\cN:=\pr_*(N_{\cY}\otimes\omega_{\cY/\SC})$ is a line bundle on $\SC$, $I_{\cY}^2\subset I_{\cYt}$, and $I_{\cY}/I_{\cYt}\cong\omega_{\cY/\SC}\otimes\pr^*\cN^{-1}$.
\item $\cCt:=\cYt\cap(H\times\SC)$ is flat over $\SC$ and is a double structure on $C\times\SC$ in $H\times\SC$ with conormal bundle $I_{C\times\SC}/I_{\cCt}\cong L_C\boxtimes\cN^{-1}$; its fibre over $[Y']$ is the ribbon $\Yt'\cap H$.
\item There is a unique morphism $\Psi\colon\SC\to\RC$ with $(\id_H\times\Psi)^{-1}(\Ct_{\RC})=\cCt$. It satisfies $\Psi^*\OO_{\RC}(-1)\cong\cN$ and, on closed points, $\Psi([Y'])=[\Yt'\cap H]$.
\end{enumerate}
\end{theorem}

\begin{proof}
(1) By Lemma~\ref{lem:SC}(2), $\SC$ is smooth at every $[Y']\in\SC$. Hence \cite[Prop.~3.7]{GP97} applies to the flat family $\pr\colon\cY\to\SC$ of smooth rational normal scrolls in $\PM\times\SC=\PP^g_{\SC}$ and gives the existence and uniqueness of $\cYt$. Its fibre over $[Y']$ is a K3 carpet on $Y'$, hence equals $\Yt'$ \cite[Prop.~1.7]{GP97}. The fact that $\cN$ is a line bundle again follows from \cite[Prop.~1.7]{GP97}. The last two claims follow from the fact that $1\in H^0(\OO_{\SC})=H^0(N_{\cY}\otimes\omega_{\cY/\SC}\otimes\pr^*\cN^{-1})$ is a nowhere vanishing section, and hence gives a surjection $N_{\cY}^{\vee}\to\omega_{\cY/\SC}\otimes\pr^*\cN^{-1}$, and $I_{\cYt}$ is the kernel of $I_{\cY}\to N_{\cY}^{\vee}\to\omega_{\cY/\SC}\otimes\pr^*\cN^{-1}$.

(2) Let $h$ be the equation of $H$. For every $[Y']\in\SC$, $h$ is a nonzerodivisor on $\OO_{Y'}$ and on $\OO_{\Yt'}$ and $Y'$ is smooth, $\Yt'$ is Cohen--Macaulay of pure dimension $2$ with support $Y'\not\subset H$. As $\cY$ and $\cYt$ are flat over $\SC$, the local criterion of flatness \cite[Cor.\ to Thm.~22.5]{Mat89} shows that $h$ is a nonzerodivisor on $\OO_{\cY}$ and on $\OO_{\cYt}$, and that $\OO_{\cCt}=\OO_{\cYt}/h\OO_{\cYt}$ is flat over $\SC$. We have $I_{C\times\SC} = \big(I_{\cY}+(h)\big)$ and $I_{\cCt} = \big(I_{\cYt}+(h)\big)$. Further, $I_{\cY}\cap(h)=hI_{\cY}$. So we have $I_{C\times\SC}^2\subset I_{\cY}^2+(h)\subset I_{\cCt}$ and
\[
I_{C\times\SC}/I_{\cCt}=\big(I_{\cY}+(h)\big)/\big(I_{\cYt}+(h)\big)=I_{\cY}/(I_{\cYt}+hI_{\cY})=(I_{\cY}/I_{\cYt})|_{C\times\SC}\cong\omega_{\cY/\SC}|_{C\times\SC}\otimes\pr^*\cN^{-1},
\]
Finally, $C\times\SC\subset\cY$ is the divisor of $h\in H^0(\OO_{\cY}(1))$, so by adjunction $\omega_{\cY/\SC}|_{C\times\SC}\cong\omega_{C\times\SC/\SC}\otimes\OO_{\cY}(-1)|_{C\times\SC}=(K_C\otimes\OO_C(-1))\boxtimes\OO_{\SC}=L_C\boxtimes\OO_{\SC}$. The fibre of $\cCt$ over $[Y']$ is $\Yt'\cap H$ because base change commutes with fibre products.

(3) By (2), $\cCt$ is defined by the surjection
\[
q\colon N_C^{\vee}\boxtimes\OO_{\SC}=I_{C\times\SC}/I_{C\times\SC}^2\longrightarrow I_{C\times\SC}/I_{\cCt}\cong L_C\boxtimes\cN^{-1},
\]
that is, $I_{\cCt}=\ker\big(I_{C\times\SC}\to N_C^{\vee}\boxtimes\OO_{\SC}\xrightarrow{\ q\ }L_C\boxtimes\cN^{-1}\big)$.
By the projection formula $q\in H^0(\SC,V\otimes\cN^{-1})$, and $q$ vanishes nowhere on $\SC$, its restriction to each $C\times\{[Y']\}$ being a surjection $N_C^{\vee}\to L_C$. Thus $q$ is a line subbundle $\cN\hookrightarrow V\otimes\OO_{\SC}$, i.e.\ a morphism $\Psi\colon\SC\to\PP(V)=\RC$ with $\Psi^*\OO_{\RC}(-1)=\cN$ and $(\id_C\times\Psi)^*\tau=q$ where $\tau$ is as in \eqref{universal ribbon}. Since $L_C\boxtimes\OO_{\RC}(1)$ is a line bundle, we have that the pullback along $\id_H\times\Psi$  of $$0\to I_{\Ct_{\RC}}\to I_{C\times\RC}\xrightarrow{\tau} L_C\boxtimes\OO_{\RC}(1)\to0$$ is the exact sequence $$0\to(\id_H\times\Psi)^*I_{\Ct_{\RC}}\to I_{C\times\SC}\xrightarrow{q}L_C\boxtimes\cN^{-1}\to0$$ and hence $(\id_H\times\Psi)^{-1}(\Ct_{\RC})=\cCt$. To see the uniqueness, note that $(\id_H\times\Psi)^{-1}(\Ct_{\RC})$ determines the surjection $N_C^{\vee}\boxtimes\OO_{\SC}\to I_{C\times\SC}/I_{\cCt}$ up to an automorphism of the target, hence the line subbundle $\Psi^*\OO_{\RC}(-1)\subset V\otimes\OO_{\SC}$, hence $\Psi$. On closed points, $\Psi([Y'])$ is the class of the fibre of $\cCt$ over $[Y']$, i.e.\ $[\Yt'\cap H]$ by (2).
\end{proof}

In the following three lemmas, we describe the map $d\Psi: H^0(N_Y(-H)) \to H^0(N_C \otimes L_C)/ks_C$, which is the differential of the map $\Psi$. An element of $H^0(N_Y(-H))$ consists of a flag $C_D = Y_{\varepsilon} \cap H_D \subset Y_{\varepsilon} \subset \mathbb{P}_D^M$ where $Y_{\varepsilon} \subset \mathbb{P}_D^M$ is an embedded first order deformation $Y \subset \mathbb{P}^M$. In the following lemma we prove a slight generalization of \cite[Prop.~3.7]{GP97}, to show that given a first order deformation $Y_{\varepsilon} \subset \mathbb{P}_D^M$ of the scroll $Y \subset \mathbb{P}^M$ there exists a unique embedded double structure $Y_{\varepsilon} \subset \widetilde Y_{\varepsilon} \subset \mathbb{P}_D^M$  that restricts on the central fiber to the unique embedded $K3$ carpet structure $Y \subset \widetilde Y \subset \mathbb{P}^M$ .

\begin{lemma}\label{lem:Phi}
Let $Y_\varepsilon\subset\PM_\Dd$ be a first-order deformation of $Y$. There is a unique closed subscheme $\Yt_\varepsilon\subset\PM_\Dd$, flat over $\Dd$, with $\Yt_\varepsilon\times_\Dd\{0\}=\Yt$, $Y_\varepsilon\subset\Yt_\varepsilon$, $I_{Y_\varepsilon}^2\subset I_{\Yt_\varepsilon}$ and $I_{Y_\varepsilon}/I_{\Yt_\varepsilon}$ invertible on $Y_\varepsilon$. Further in this case, $I_{Y_\varepsilon}/I_{\Yt_\varepsilon} \cong \omega_{Y_{\varepsilon}/D}$.  We write $\Phi(Y_\varepsilon):=\Yt_\varepsilon$ and $d\Phi\colon H^0(\NY)\to H^0(\NYt)$ for the resulting map. Moreover:
\begin{enumerate}
\item the image of $d\Phi(x)$ under $H^0(\NYt) = \operatorname{Hom}(I_{\widetilde Y}, \mathcal O_{\widetilde Y})\to H^0(\NYt\otimes\OO_Y) = \operatorname{Hom}(I_{\widetilde Y}, \mathcal O_{Y})$ is $\iota(x):=x|_{I_{\Yt}}$, the image of $x$ under the composite $H^0(\NY) = \operatorname{Hom}(I_{Y}, \mathcal O_{Y})\to H^0(\Hc)\to H^0(\NYt\otimes\OO_Y) = \operatorname{Hom}(I_{\widetilde Y}, \mathcal O_{Y})$ of the maps in \eqref{eq:3.3} and \eqref{eq:3.2};
\item for $\gamma_\varepsilon\in PGL_{M+1}(\Dd)$ one has $\Phi(\gamma_\varepsilon Y)=\gamma_\varepsilon\Yt$.
\end{enumerate}
\end{lemma}

\begin{proof}
Such a $\Yt_\varepsilon$ is the kernel of a surjection $s_\varepsilon\colon N_{Y_\varepsilon/\PM_\Dd}^{\vee}\to L_\varepsilon$ onto an invertible sheaf with $L_\varepsilon|_Y\cong K_Y$. Since $H^1(\mathcal{O}_Y) = 0$, we have that $\operatorname{Pic}(Y_\varepsilon)\to\operatorname{Pic}(Y)$ is injective. Hence  $L_\varepsilon\cong\omega_{Y_\varepsilon/\Dd}$, and $s_\varepsilon$ is determined by $\Yt_\varepsilon$ up to $H^0(\OO_{Y_\varepsilon}^{\times})=(k[\varepsilon]/\varepsilon^2)^{\times}$. Conversely every such surjection lifting $s$ defines such a $\Yt_\varepsilon$. Put $\mathcal{F}:=N_{Y_\varepsilon/\PM_\Dd}\otimes\omega_{Y_\varepsilon/\Dd}$, a locally free sheaf on $Y_\varepsilon$. Tensoring the exact sequence 
$$0 \to \varepsilon \mathcal{O}_Y \to \mathcal{O}_{Y_{\varepsilon}} \to \mathcal{O}_Y \to 0$$
by $\mathcal F$, we have $$0\to\varepsilon\mathcal{F}\to\mathcal{F}\to\NY\otimes K_Y\to 0$$ 
and $\varepsilon\mathcal{F}\cong\NY\otimes K_Y$. By Lemma~\ref{lem:pushforward}, $H^0(\NY\otimes K_Y)=k\,s$,  $H^1(\NY\otimes K_Y)=0$. Therefore, $H^0(\mathcal{F})\to H^0(\NY\otimes K_Y)=k\,s$ is surjective with kernel $\varepsilon H^0(\mathcal{F})=\varepsilon k\,s$. So lifts of $s$ exist, are surjective by Nakayama, and form the set $\{(1+\varepsilon\lambda)s_\varepsilon\}$. But every such surjection has the same kernel and therefore we have a unique double structure $\widetilde{Y}_{\varepsilon}$.

(1) With $\eta:=d\Phi(x)$, $I_{\Yt_\varepsilon}\subset I_{Y_\varepsilon}$ means that for $f\in I_{\Yt}$ and any $g$ with $g\equiv\eta(f)$ mod $I_{\Yt}$ one has $g\equiv x(f)$ mod $I_Y$, i.e.\ $\eta(f)\equiv x(f)$ mod $I_Y$ for all $f\in I_{\Yt}$. (2) follows from uniqueness.
\end{proof}

\begin{lemma}\label{def:Psi}
Let $x\in H^0(\NY(-H))=T_{[Y]}\SC$, let $Y_\varepsilon\supset C_\Dd$ be the corresponding deformation (Lemma~\ref{lem:SC}(1)), $\Yt_\varepsilon:=\Phi(Y_\varepsilon)$, and
\[
\Ct_\varepsilon:=\Yt_\varepsilon\cap H_\Dd\subset H_\Dd
\]
(scheme-theoretic intersection). Then:
\begin{enumerate}
\item[(i)] $\Ct_\varepsilon$ is flat over $\Dd$ with closed fibre $\Ct$. 
\item[(ii)] $\Ct_\varepsilon$ is a double structure on $C_\Dd$ with conormal bundle $L_C\otimes\OO_\Dd$. 
\end{enumerate}
Hence in particular, by Lemma~\ref{lem:RC}(2), $\Ct_\varepsilon=\ker(s_C+\varepsilon t)$ for a $t\in H^0(N_C\otimes L_C)$ well defined modulo $k\,s_C$.
\end{lemma}

\begin{proof}
To see (i) note that $h$ is a nonzerodivisor on $\OO_{\Yt}$, since $\Yt$ is Cohen--Macaulay of pure dimension $2$ and $h$ is a nonzerodivisor on $\OO_{Y}$. Now the conclusion follows from local criterion of flatness. \\
To see (ii), note that since $\widetilde Y_{\varepsilon}$ is double structure over $Y_{\varepsilon}$ with conormal bundle $\omega_{Y_{\varepsilon}/D}$ by Lemma~\ref{lem:Phi}, $\widetilde C_{\varepsilon} = \widetilde Y_{\varepsilon} \cap H_D$ is a double structure on $Y_{\varepsilon} \cap H_D = C_D$ with conormal bundle $\omega_{Y_{\varepsilon}/D}|_{C_D} = L_C \otimes \mathcal O_D$. 
\end{proof}

\begin{lemma}\label{lem:dPsi}
The differential
\[
d\Psi\colon T_{[Y]}\SC=H^0(\NY(-H))\longrightarrow H^0(N_C\otimes L_C)/k\,s_C=T_{[\Ct]}\RC
\]
of the map $\Psi$ of Theorem~\ref{thm:Psi} is given by
\[
d\Psi(Y_\varepsilon\supset C_\Dd)=\Ct_\varepsilon=\Yt_\varepsilon\cap H_\Dd .
\]
In other words, if $x\in H^0(\NY(-H))$ and $\Ct_\varepsilon=\ker(s_C+\varepsilon t)$ for a $t\in H^0(N_C\otimes L_C)$, then $d\Psi(x)=[t]$, where $[t]$ denotes the class of $t$ modulo $k\,s_C$.
\end{lemma}

\begin{proof}
 Regard $Y_\varepsilon \supset C_D$ as a morphism $\Dd\to\SC$. Then $d\Psi(Y_\varepsilon \supset C_D)$ is given by the composition $D \to \mathcal{S}_C \xrightarrow{\Psi} \mathcal{R}_C$. By Theorem~\ref{thm:Psi}(1) and flatness, $\cYt\times_{\SC}\Dd$ is flat over $\Dd$ and is a double structure on $Y_\varepsilon=\cY\times_{\SC}\Dd$ such that its restriction to the central fiber is the $K3$ double structure $\widetilde Y $ on $Y$. So
 $\cYt\times_{\SC}\Dd=\Yt_\varepsilon$ by the uniqueness in Lemma~\ref{lem:Phi}. Hence $\Ct_\varepsilon=\cCt\times_{\SC}\Dd$ is the pullback of $\Ct_{\RC}$ along $D \to \mathcal{S}_C \xrightarrow{\Psi} \mathcal{R}_C$ by Theorem~\ref{thm:Psi}(3). Therefore $d\Psi(Y_\varepsilon\supset C_\Dd)=\Ct_\varepsilon=\Yt_\varepsilon\cap H_\Dd$. 
\end{proof}

\begin{proposition}\label{prop:a=dPsi}
$\bar\delta^{-1} \circ a= d\Psi$. Equivalently, the map \eqref{eq:themap} is $d\Psi$.
\end{proposition}

\begin{proof}
Throughout we regard $H^1(\NY(-H+K_Y))$ as a subspace of $H^1(\NYt(-H+K_Y))$ via the injective map $j$ of the bottom row of \eqref{diag:a}. The map $j$ is induced by the sheaf map $\NY\to\NYt\otimes\OO_Y$, which restricts a homomorphism from $I_Y$ to $I_{\Yt}$ (the composite of \eqref{eq:3.3} and \eqref{eq:3.2}), and so is the injection $\iota\colon H^0(\NY(-H))\hookrightarrow H^0(\NYt(-H))$ of \eqref{eq:3.60}, \eqref{eq:3.61}. The maps in play are
\begin{equation}\label{diag:overview}
\begin{tikzcd}[column sep=large, row sep=large]
&& H^0(\NYt(-\Ht))\arrow[d]\\
0\arrow[r] & H^0(\NY(-H))\arrow[r,"\iota"]\arrow[d,dashed,"x\mapsto\eta|_{\Ct}"]\arrow[dr,"a"' near end]\arrow[dd,bend right=60,"d\Psi"'] & H^0(\NYt(-H))\arrow[d,"F"]\\
& H^0(N_C\otimes L_C)\arrow[r,"\delta"]\arrow[d,two heads] & H^1(\NYt(-H+K_Y))\\
& H^0(N_C\otimes L_C)/k\,s_C=T_{[\Ct]}\RC\arrow[ur,"\bar\delta"',"\sim"] &
\end{tikzcd}
\end{equation}
Here $F$ is the connecting map of \eqref{eq:3.65}, i.e.\ of \eqref{eq:3.1} twisted by $-\Ht$, and by definition of $a$ (diagram \eqref{diag:a}) $a=F\circ\iota$. The row through $\delta$ is the cohomology of the restriction sequence $$0\to\NY\otimes K_Y(-C)\to\NY\otimes K_Y\to N_C\otimes L_C\to0$$ of Lemma~\ref{lem:delta}. $\delta$ is its connecting map, its kernel $k\,s_C$ is the image of $H^0(\NY\otimes K_Y)=k\,s$, and $\bar\delta$ is the induced isomorphism. The dashed arrow is explained below; it is well defined modulo $k\,s_C=\ker\delta$.

Pick $x\in H^0(\NY(-H))$ and set $\xi:=\iota(x)$. We need to show $\bar\delta(d\Psi(x))=F(\xi)$.

%\emph{Computation of $d\Psi(x)$.} 
The vector $x$ is a first-order deformation $Y_\varepsilon\subset\PM_\Dd$ of $Y$ containing $C_\Dd$, and $Y_\varepsilon\cap H_\Dd=C_\Dd$ (Lemma~\ref{lem:SC}(1), Lemma~\ref{def:Psi}(ii)). By Lemma~\ref{lem:Phi} (the version of \cite[Prop.~3.7]{GP97} over $\Dd$) there is a K3 carpet $\Yt_\varepsilon=\Phi(Y_\varepsilon)\subset\PM_\Dd$ on $Y_\varepsilon$; since $\Yt_\varepsilon\times_\Dd\operatorname{Spec}k=\Yt$, it is an embedded first-order deformation of $\Yt\subset\PM$. Call its class $\eta\in H^0(\NYt)$. Then $\Ct_\varepsilon:=\Yt_\varepsilon\cap H_\Dd$ is a ribbon on $C_\Dd$ with conormal bundle $L_C\otimes\OO_\Dd$, hence a $\Dd$-point of $\RC$ at $[\Ct]$, and by definition (Lemma~\ref{def:Psi})
\[
d\Psi(x)=[\Ct_\varepsilon]\in T_{[\Ct]}\RC=H^0(N_C\otimes L_C)/k\,s_C .
\]

%\emph{Where $\eta|_{\Ct}$ lives.} 
Put $N_{\Ct}:=N_{\Ct/H}=\sHom(I_{\Ct},\OO_{\Ct})$. Composing $\eta\colon I_{\Yt}\to\OO_{\Yt}$ with $\OO_{\Yt}\to\OO_{\Ct}$ and using $I_{\Yt}\otimes\OO_H=I_{\Ct/H}$ gives $\eta|_{\Ct}\in H^0(N_{\Ct})$. The class of $\Ct_\varepsilon$ as a first-order deformation of $\Ct$ inside $H$ (reduce the generators $f+\varepsilon\tilde\eta(f)$ of $I_{\Yt_\varepsilon}$ modulo $h$). For the ribbon $\Ct\subset H$, whose conormal bundle is $L_C=I_C/I_{\Ct}$, we have the analogues of \eqref{eq:3.1}--\eqref{eq:3.3} \cite[Lemma~4.1]{GGP08}:
\begin{equation}\label{eq:ribbonseq}
\begin{gathered}
0\to N_{\Ct}\otimes L_C\to N_{\Ct}\to N_{\Ct}\otimes\OO_C\to0,\qquad
0\to\Hc_C\to N_{\Ct}\otimes\OO_C\to L_C^{-2}\to0,\\
0\to L_C^{-1}\xrightarrow{\ s_C^{\vee}\ }N_C\to\Hc_C\to0,\qquad\text{where }\Hc_C:=\sHom(I_{\Ct}/I_C^2,\OO_C).
\end{gathered}
\end{equation}

\emph{Claim.}
\begin{enumerate}
\item[(i)] $\eta|_{\Ct}\in H^0(N_{\Ct}\otimes L_C)\subset H^0(N_{\Ct})$.
\item[(ii)] Restriction of homomorphisms from $I_C$ to $I_{\Ct}$ induces an injection $H^0(N_C\otimes L_C)/k\,s_C\hookrightarrow H^0(N_{\Ct}\otimes L_C)$, and $\eta|_{\Ct}$ lies in its image.
\item[(iii)] Under this injection, $\eta|_{\Ct}=-d\Psi(x)$.
\end{enumerate}

\emph{Proof of the Claim.} (i) The image of $\eta|_{\Ct}$ in $H^0(N_{\Ct}\otimes\OO_C)=\Hom(I_{\Ct},\OO_C)$ is $\eta|_C\colon f\mapsto\eta(f)|_C$. By Lemma~\ref{lem:Phi}(1), $\eta(f)\equiv x(f)$ mod $I_Y$; so $\eta|_C$ is the restriction to $I_{\Ct}$ of the class $f\mapsto x(f)|_C$ of $Y_\varepsilon\cap H_\Dd$ as a deformation of $C$ in $H$, which is $0$ because $Y_\varepsilon\cap H_\Dd=C_\Dd$. Hence $\eta|_C=0$, and by the first sequence of \eqref{eq:ribbonseq}, $\eta|_{\Ct}\in H^0(N_{\Ct}\otimes L_C)=\Hom(I_{\Ct},L_C)$.

(ii) Tensoring the third sequence of \eqref{eq:ribbonseq} with $L_C$ gives $0\to\OO_C\xrightarrow{s_C}N_C\otimes L_C\to\Hc_C\otimes L_C\to0$, whence $H^0(N_C\otimes L_C)/k\,s_C\cong H^0(\Hc_C\otimes L_C)$ as $H^1(\OO_C)=0$; tensoring the second with $L_C$ gives $0\to\Hc_C\otimes L_C\to N_{\Ct}\otimes L_C\to L_C^{-1}\to0$, whence $H^0(\Hc_C\otimes L_C)\hookrightarrow H^0(N_{\Ct}\otimes L_C)$. The composite is $t\mapsto t|_{I_{\Ct}}$. Now $\Ct_\varepsilon$ is a double structure on $C_\Dd$, i.e.\ $I_{C_\Dd}^2\subset I_{\Ct_\varepsilon}$ (Lemma~\ref{def:Psi}(ii)); in terms of the class $\eta|_{\Ct}$ this says exactly that $\eta|_{\Ct}\in\Hom(I_{\Ct},L_C)$ vanishes on $I_C^2$, i.e.\ $\eta|_{\Ct}\in\Hom(I_{\Ct}/I_C^2,L_C)=H^0(\Hc_C\otimes L_C)$.

(iii) Write $d\Psi(x)=[t]$, i.e.\ $\Ct_\varepsilon=\ker(s_C+\varepsilon t)$ (Lemma~\ref{def:Psi}). For $f\in I_{\Ct}$ the element $f+\varepsilon g$ lies in $I_{\Ct_\varepsilon}$ if and only if $s_C(g)=-t(f)$; so the class of $\Ct_\varepsilon$ is $f\mapsto-t(f)\in L_C=I_C/I_{\Ct}$, that is, $\eta|_{\Ct}=-t|_{I_{\Ct}}$. This proves the Claim.

%\emph{The two connecting maps.} 
Now let $\tilde\delta\colon H^0(N_{\Ct}\otimes L_C)\to H^1(\NYt(-H+K_Y))$ be the connecting map of the restriction sequence
\[
0\to\NYt(-H+K_Y)\to\NYt(K_Y)\to\NYt(K_Y)|_C=N_{\Ct}\otimes L_C\to0
\]
of $\NYt(K_Y)=\sHom(I_{\Yt},K_Y)$, $K_Y=I_Y/I_{\Yt}$. Restriction of homomorphisms from $I_Y$ to $I_{\Yt}$ maps the restriction sequence of $\NY\otimes K_Y$ to this one, inducing $j$ on $H^1$ of the first terms and the injection of (ii) on $H^0$ of the third terms; so $\tilde\delta$ restricts to $\delta$ on $H^0(N_C\otimes L_C)/k\,s_C$, and by (iii)
\[
\tilde\delta(\eta|_{\Ct})=-\bar\delta(d\Psi(x)) .
\]
Therefore what we need to show is the following.

\emph{Lemma.} $\tilde\delta(\eta|_{\Ct})=-F(\xi)$.

\emph{Proof of the Lemma.} Multiplication by $h$ followed by restriction to $\Ct$ gives a commutative diagram of sheaves on $\Yt$ with exact rows and columns
\begin{equation}\label{diag:3x3}
\begin{tikzcd}[row sep=small, column sep=small]
& 0\arrow[d] & 0\arrow[d] & 0\arrow[d] & \\
0\arrow[r] & \NYt(-H+K_Y)\arrow[r]\arrow[d,"h"'] & \NYt(-\Ht)\arrow[r]\arrow[d,"h"'] & \NYt(-H)\arrow[r]\arrow[d,"h"] & 0\\
0\arrow[r] & \NYt(K_Y)\arrow[r]\arrow[d] & \NYt\arrow[r]\arrow[d] & \NYt\otimes\OO_Y\arrow[r]\arrow[d] & 0\\
0\arrow[r] & N_{\Ct}\otimes L_C\arrow[r]\arrow[d] & N_{\Ct}\arrow[r]\arrow[d] & N_{\Ct}\otimes\OO_C\arrow[r]\arrow[d] & 0\\
& 0 & 0 & 0 &
\end{tikzcd}
\end{equation}
whose rows are \eqref{eq:3.1} twisted by $-\Ht$, \eqref{eq:3.1} itself, and the first sequence of \eqref{eq:ribbonseq}; the columns are exact because $h$ is a nonzerodivisor on $\OO_{\Yt}$ and on $\OO_Y$. $F$ is the connecting map of the top row and $\tilde\delta$ that of the left column. Lift $\xi$ locally to sections $\xi_i$ of $\NYt(-\Ht)$ over an affine cover $\{U_i\}$; then $F(\xi)=[\xi_i-\xi_j]$. We make two observations:
\begin{enumerate}
\item[(a)] View the $\xi_i$ as local sections of $\NYt$ (middle column). Then $\eta|_{U_i}-\xi_i$ restricts to zero on $Y$, since $\eta$ and $\xi_i$ both reduce to $\xi=x|_{I_{\Yt}}$ (Lemma~\ref{lem:Phi}(1)). Hence $\mu_i:=\eta|_{U_i}-\xi_i$ are local sections of $\NYt(K_Y)$.
\item[(b)] The $\xi_i$ vanish on $\Ct$, being multiples of $h$; so $\mu_i|_{\Ct}=\eta|_{\Ct\cap U_i}$, i.e.\ the $\mu_i$ are local lifts of $\eta|_{\Ct}$ to $\NYt(K_Y)$.
\end{enumerate}
Consequently $\tilde\delta(\eta|_{\Ct})=[\mu_i-\mu_j]=[\xi_j-\xi_i]=-F(\xi)$, which proves the Lemma.

By the Lemma, $\bar\delta(d\Psi(x))=-\tilde\delta(\eta|_{\Ct})=F(\xi)=a(x)$.
\end{proof}

\subsubsection{Construction of local sections of \texorpdfstring{$\Psi$}{Psi} and surjectivity of the differential \texorpdfstring{$d\Psi$}{dPsi}}\label{subsec:sections}

In this subsubsection, we construct the local sections of $\Psi$ and show that the differential $d\Psi$ is surjective.
Put $Q:=K_C\otimes L_C^{-1}=\OO_C(C)=\OO_{\PP^1}(g-1)$. For $Y'\in\SC$ let $\pi'\colon Y'\to C$ be the ruling for which $C$ is a section and $E_{Y'}:=\pi'_*\OO_{Y'}(C)$. Pushing forward we have $0\to\OO_{Y'}\to\OO_{Y'}(C)\to\OO_C(C)\to0$ since $R^1\pi'_*\OO_{Y'}=0$. This gives the extension
\begin{equation}\label{eq:EY}
0\to\OO_C\to E_{Y'}\to Q\to0,
\end{equation}
and $Y'=\PP(E_{Y'})$ with $\OO_{\PP(E_{Y'})}(1)=\OO_{Y'}(C)$ and $C$ the section defined by $E_{Y'}\to Q$. Since $Y'$ is smooth, $E_{Y'}\cong\OO(a_1')\oplus\OO(a_2')$ with $a_i'\ge1$, so \eqref{eq:EY} is not split. Put $V_{Y'}:=E_{Y'}\otimes L_C$, an extension $0\to L_C\to V_{Y'}\to K_C\to0$ with class $[V_{Y'}]\in\Ext^1(K_C,L_C)=\Ext^1(\Omega_C,L_C)$, $[V_{Y'}]\ne0$. We will also denote the extension class of \eqref{eq:EY} by $[E_{Y'}]$.

The two spaces $\RC$ and $\PP(\Ext^1(\Omega_C,L_C))$ are identified by the classical dictionary of Bayer--Eisenbud:

\begin{lemma}\label{lem:beta}
Let $\varepsilon_C\in\Ext^1(\Omega_C,N_C^{\vee})$ be the class of the conormal sequence $0\to N_C^{\vee}\to\Omega_H|_C\to\Omega_C\to0$, and let
\[
\beta\colon H^0(N_C\otimes L_C)=\Hom(N_C^{\vee},L_C)\longrightarrow\Ext^1(\Omega_C,L_C),\qquad t\longmapsto t_*(\varepsilon_C)
\]
(push-out along $t$). Then $\beta(t)$ is the class of the restricted cotangent sequence $0\to L_C\to\Omega_{\Ct_t}|_C\to\Omega_C\to0$ of the double structure $\Ct_t$, and $\beta$ is an isomorphism of $(g-2)$-dimensional vector spaces.
\end{lemma}

\begin{proof}
Tensoring the conormal sequence of $\Ct_t\subset H$ with $\OO_C$ gives $I_{\Ct_t}/I_CI_{\Ct_t}\to\Omega_H|_C\to\Omega_{\Ct_t}|_C\to 0$, and the image of the first map is the image of $\ker t=I_{\Ct_t}/I_C^2\subset N_C^{\vee}$ in $\Omega_H|_C$. Hence $\Omega_{\Ct_t}|_C=\Omega_H|_C/\ker t$, and $0\to N_C^{\vee}/\ker t=L_C\to\Omega_H|_C/\ker t\to\Omega_C\to0$ is the push-out of the conormal sequence along $t$; this is the first claim (cf.\ \cite[Theorem~1.2]{BE95}). For the second, apply $\Hom(-,L_C)$ to the conormal sequence: $\ker\beta$ is the image of $\Hom(\Omega_H|_C,L_C)=H^0(T_H|_C\otimes L_C)$. From the Euler sequence of $H\cong\PP^{g-1}$ restricted to $C$ and tensored with $L_C$,
\[
0\to L_C\to L_C(1)^{\oplus g}\to T_H|_C\otimes L_C\to 0,\qquad L_C(1)=\OO_{\PP^1}(-2),
\]
we get $H^0(T_H|_C\otimes L_C)=\ker\big(H^1(L_C)\to H^1(L_C(1))\otimes H^0(\OO_H(1))^{\vee}\big)$. This map is adjoint to the multiplication pairing $H^1(\OO_{\PP^1}(-g-1))\otimes H^0(\OO_{\PP^1}(g-1))\to H^1(\OO_{\PP^1}(-2))=k$, which is perfect by Serre duality; hence it is injective and $H^0(T_H|_C\otimes L_C)=0$. So $\beta$ is injective, and both spaces have dimension $g-2$.
\end{proof}

In the following theorem, given a point in $\mathcal{R}_C$ we construct a map $\Gamma_U$ defined around a neighbourhood $U$ of the point. We will show in Proposition~\ref{prop:BEsection}, that $\Gamma_U$ is a section of $\Psi$ over $U$.

\begin{theorem}\label{thm:Gamma}
Let $\pr\colon C\times\RC\to\RC$ be the projection.
\begin{enumerate}
\item The restricted cotangent sequence of the universal ribbon $\Ct_{\RC}\subset H\times\RC$,
\begin{equation}\label{eq:relcot}
0\to L_C\boxtimes\OO_{\RC}(1)\to\Omega_{\Ct_{\RC}/\RC}|_{C\times\RC}\to K_C\boxtimes\OO_{\RC}\to0,
\end{equation}
is exact, and $\mathcal{V}:=\Omega_{\Ct_{\RC}/\RC}|_{C\times\RC}\otimes(L_C^{-1}\boxtimes\OO_{\RC}(-1))$ sits in an extension
\begin{equation}\label{eq:VR}
0\to\OO_{C\times\RC}\to\mathcal{V}\to Q\boxtimes\OO_{\RC}(-1)\to0
\end{equation}
whose fibre $\mathcal{V}_u$ over $u=[t]\in\RC$ is isomorphic to $\Omega_{\Ct_t}|_C\otimes L_C^{-1}$: a non-split extension of $Q$ by $\OO_C$ with $[\mathcal{V}_u\otimes L_C]\in k^{\times}\beta(t)$.
\item Every point of $\RC$ has an open neighbourhood $U$ with a morphism $\Gamma_U\colon U\to\SC$ such that, for every $u\in U$, $\Gamma_U(u)$ is $\PP(\mathcal{V}_u)$ embedded by $|\OO_{\PP(\mathcal{V}_u)}(1)|$, with $C$ the section defined by $\mathcal{V}_u\to Q$. In particular $[E_{\Gamma_U(u)}]=[\mathcal{V}_u]$ as extensions of $Q$ by $\OO_C$, and $[V_{\Gamma_U(u)}]\in k^{\times}\beta(t)$ for $u=[t]$.
\end{enumerate}
\end{theorem}

\begin{proof}
(1) The sequence \eqref{eq:relcot} is the conormal sequence of $C\times\RC\subset\Ct_{\RC}$ relative to $\RC$: its first term is $I_{C\times\RC}/I_{\Ct_{\RC}}=L_C\boxtimes\OO_{\RC}(1)$ and $\Omega_{C\times\RC/\RC}=K_C\boxtimes\OO_{\RC}$. It is exact on the left because $\Ct_{\RC}$ is locally isomorphic to $(C\times\RC)\times\Dd$ (ribbons are locally split, \cite[\S1]{BE95}; the same local computation applies to $\Ct_{\RC}$, which is defined by the surjection $\tau$). Its restriction to $C\times\{[t]\}$ is the restricted cotangent sequence $0\to L_C\to\Omega_{\Ct_t}|_C\to K_C\to0$ of $\Ct_t$, of class $\beta(t)\ne0$ by Lemma~\ref{lem:beta}. Tensoring with $L_C^{-1}\boxtimes\OO_{\RC}(-1)$ gives \eqref{eq:VR} and the description of its fibres.

(2) The surjection $\mathcal{V}\to Q\boxtimes\OO_{\RC}(-1)$ gives a section $i\colon C\times\RC\to\PP(\mathcal{V})$ of $p\colon\PP(\mathcal{V})\to C\times\RC$ with $i^*\OO_{\PP(\mathcal{V})}(1)=Q\boxtimes\OO_{\RC}(-1)$; since the kernel of this surjection is $\OO_{C\times\RC}$, $i(C\times\RC)$ is the divisor of a section of $\OO_{\PP(\mathcal{V})}(1)$. Let $\mathcal{W}:=\pr_*p_*\OO_{\PP(\mathcal{V})}(1)=\pr_*\mathcal{V}$. For $u\in\RC$, $\mathcal{W}_u=H^0(\mathcal{V}_u)$ with $\mathcal{V}_u\cong\OO(a_1)\oplus\OO(a_2)$, $a_1+a_2=g-1$, $a_i\ge1$ (an extension of $\OO(g-1)$ by $\OO$ with $a_1=0$ is split); so $h^0(\mathcal{V}_u)=g+1=M+1$, $h^1(\mathcal{V}_u)=0$, and $\mathcal{W}$ is a vector bundle of rank $M+1$. Embed $\PP(\mathcal{V})$ in $\PP(\mathcal{W})$ over $\RC$ by $|\OO_{\PP(\mathcal{V})}(1)|$. Pushing forward \eqref{eq:VR},
\[
0\to\OO_{\RC}\to\mathcal{W}\to H^0(Q)\otimes\OO_{\RC}(-1)\to0,
\]
and $i(C\times\RC)$ lies in $\PP\big(H^0(Q)\otimes\OO_{\RC}(-1)\big)=\PP(H^0(Q))\times\RC=H\times\RC$, where it is $C\times\RC$, because $i(C\times\{u\})$ is embedded by $|i^*\OO(1)|=|\OO_C(1)|$. Given $u\in\RC$, let $U\ni u$ be an open set on which $\OO_{\RC}(-1)$ is trivial, e.g.\ the complement of a hyperplane, and choose a splitting of the sequence (it exists, as $\Ext^1(H^0(Q)\otimes\OO_{\RC}(-1),\OO_{\RC})=H^0(Q)^{\vee}\otimes H^1(\OO_{\RC}(1))=0$). This gives $\mathcal{W}|_U\cong(k\oplus H^0(Q))\otimes\OO_U=H^0(\OO_{\PM}(1))\otimes\OO_U$ compatible with the quotients onto $H^0(Q)=H^0(\OO_H(1))$, hence $\PP(\mathcal{W})|_U\cong\PM\times U$ carrying $H\times U$ to $H\times U$ and $C\times U$ to $C\times U$. The composite $\PP(\mathcal{V})|_U\hookrightarrow\PP(\mathcal{W})|_U\cong\PM\times U$ is a flat family of smooth scrolls containing $C\times U$, i.e.\ a morphism $\Gamma_U\colon U\to\SC$. By construction $\Gamma_U(u)=\PP(\mathcal{V}_u)$ embedded by $|\OO_{\PP(\mathcal{V}_u)}(1)|$ with $C=i(C\times\{u\})$; as $\OO_{\PP(\mathcal{V}_u)}(1)=\OO(C)$, $E_{\Gamma_U(u)}=\pi_*\OO_{\Gamma_U(u)}(C)=\mathcal{V}_u$, and \eqref{eq:EY} for $\Gamma_U(u)$ is the fibre of \eqref{eq:VR}. The last assertion follows from (1).
\end{proof}

%\subsubsection*{$\GH$-invariance}
In the following lemma we show the $G_H$ invariance of $\Psi$, where $G_H$ is the set of automorphisms which induces identity on $H$.

\begin{lemma}\label{lem:GH}
\begin{enumerate}
\item $\GH$ acts on $\SC$, and $\Psi\circ\gamma=\Psi$ for every $\gamma\in\GH$.
\item If $Y',Y''\in\SC$ and the classes $[V_{Y'}],[V_{Y''}]\in\Ext^1(K_C,L_C)$ are proportional, then $Y''=\gamma Y'$ for some $\gamma\in\GH$.
\end{enumerate}
\end{lemma}

\begin{proof}
(1) A projectivity $\gamma$ with $\gamma|_H=\id$ maps $Y'\in\SC$ to a smooth scroll $\gamma Y'\supset\gamma C=C$, so $\GH$ acts on $\SC$. Next, $\gamma\Yt'$ is a K3 carpet on $\gamma Y'$, hence is the K3 carpet on $\gamma Y'$ \cite[Prop.~1.7]{GP97}, and $\gamma\Yt'\cap H=\gamma(\Yt'\cap H)=\Yt'\cap H$ as subschemes of $H$ because $\gamma|_H=\id$. Thus $\Psi(\gamma Y')=\Psi(Y')$ on closed points, and $\Psi\circ\gamma=\Psi$ since $\SC$ is reduced (proof of Theorem~\ref{thm:Psi}).

(2) Proportional classes give an isomorphism $E_{Y'}\cong E_{Y''}$ inducing the identity on $Q$ (and a nonzero scalar on $\OO_C$), hence an isomorphism $(\PP(E_{Y'}),C)\cong(\PP(E_{Y''}),C)$ over $C$. So $Y'$ and $Y''$ are two embeddings of the same pair by the complete linear series $|\OO(C)|$ which agree on $C$; they differ by a $\gamma\in PGL_{M+1}(k)$ with $\gamma|_C=\id$, hence $\gamma|_H=\id$ since $C$ spans $H$, i.e.\ $\gamma\in\GH$.
\end{proof}

In the following Proposition, we show that $\Gamma_U$ is a local section of $\Psi$.

\begin{proposition}\label{prop:BEsection}
Let $\Gamma_U\colon U\to\SC$ be as in Theorem~\ref{thm:Gamma}(2).
\begin{enumerate}
\item $\Psi\circ\Gamma_U=\id_U$.
\item For every $Y'\in\SC$, $\Psi(Y')=[t]$ for the $t\in H^0(N_C\otimes L_C)$ with $\beta(t)\in k^{\times}[V_{Y'}]$, i.e., the hyperplane section $\Yt'\cap H$ of the carpet is the ribbon whose restricted cotangent extension is that of $[V_{Y'}]=[E_{Y'}\otimes L_C]$.
\item $\Gamma_U(\Psi(Y'))\in\GH\cdot Y'$ for every $Y'\in\Psi^{-1}(U)$.
\end{enumerate}
\end{proposition}

\begin{proof}
(1) Let $u=[t]\in U$ and $\Ct:=\Ct_t\subset H$. It is a non-split ribbon on $C$ with conormal bundle $L_C$ (Lemma~\ref{lem:RC}, Lemma~\ref{lem:beta}), and it is canonically embedded in $H=\PP^{g-1}$: it is nondegenerate of degree $2g-2$ and arithmetic genus $g$ in $\PP^{g-1}$, so by Riemann--Roch $\omega_{\Ct}\otimes\OO_{\Ct}(-1)$ is a line bundle of degree $0$ with a nonzero section, hence trivial (its restriction to $C$ is trivial, and $H^0(L_C)=0$). We apply Theorem~\ref{thm:extension}(2) to $\Ct\subset H$ (so $g_C=0$ in Set-up~\ref{setup:ribbon}). By Theorem~\ref{thm:Gamma}(1) the restricted cotangent sequence of $\Ct$ is the fibre of \eqref{eq:relcot} at $u$, so $V=\Omega_{\Ct}|_C=\mathcal V_u\otimes L_C$; hence the ruled surface of Set-up~\ref{setup:ribbon} is $\PP(V)=\PP(\mathcal V_u)$, with $\OO_{\PP(V)}(C)=\OO_{\PP(V)}(1)\otimes p^*L_C^{-1}=\OO_{\PP(\mathcal V_u)}(1)$ (Lemma~\ref{lem:OSC}) and $C$ the section defined by $\mathcal V_u\to Q$. By Theorem~\ref{thm:Gamma}(2), $\Gamma_U(u)\subset\PM$ is $\PP(\mathcal V_u)$ embedded by the complete linear series of $\OO_{\PP(\mathcal V_u)}(1)$, with $H$ the hyperplane such that $\Gamma_U(u)\cap H=C$ is the given curve (Lemma~\ref{lem:SC}); and $\Yt'$ is the K3 carpet on $\Gamma_U(u)$. Theorem~\ref{thm:extension}(2) gives $\Yt'\cap H=\Ct$, hence $\Psi(\Gamma_U(u))=[\Yt'\cap H]=[\Ct]=u$.

(2), (3) Let $Y'\in\SC$. Since $[E_{Y'}]$ is non-split, $[V_{Y'}]\ne0$, and since $\beta$ is an isomorphism (Lemma~\ref{lem:beta}) there is a $t\ne0$ with $\beta(t)=[V_{Y'}]$. Put $u:=[t]$ and choose $\Gamma_{U'}$ with $u\in U'$ (Theorem~\ref{thm:Gamma}(2)). Then $[V_{\Gamma_{U'}(u)}]\in k^{\times}\beta(t)=k^{\times}[V_{Y'}]$, so $Y'\in\GH\cdot\Gamma_{U'}(u)$ by Lemma~\ref{lem:GH}(2), and by Lemma~\ref{lem:GH}(1) and part (1),
\[
\Psi(Y')=\Psi(\Gamma_{U'}(u))=u=[t] .
\]
This is (2). If moreover $\Psi(Y')\in U$, then $u=\Psi(Y')\in U$, and $[V_{\Gamma_U(u)}]$ also lies in $k^{\times}\beta(t)$, so $\Gamma_U(u)\in\GH\cdot\Gamma_{U'}(u)=\GH\cdot Y'$ by Lemma~\ref{lem:GH}(2). This is (3).
\end{proof}

Finally in the following Proposition, we show that $d\Psi$ is surjective.

\begin{proposition}\label{prop:dPsi}
Let $d\Psi\colon H^0(\NY(-H))=T_{[Y]}\SC\to T_{[\Ct]}\RC=H^0(N_C\otimes L_C)/k\,s_C$ be the differential of $\Psi$ at $[Y]$ (see Lemma~\ref{def:Psi}).
\begin{enumerate}
\item $d\Psi$ is surjective.
\item $\ker d\Psi=\Im c=T_{[Y]}(\GH\cdot Y)$.
\end{enumerate}
\end{proposition}

\begin{proof}
(1) Choose an open $U\ni\Psi(Y)$ and $\Gamma_U$ as in Theorem~\ref{thm:Gamma}(2). By Proposition~\ref{prop:BEsection}(1), $\Psi\circ\Gamma_U=\id_U$; differentiating at $\Psi(Y)$,
\[
(d\Psi)_{\Gamma_U(\Psi(Y))}\circ {(d\Gamma_U)}_{\Psi(Y)}=\id\colon T_{\Psi(Y)}\RC\to T_{\Psi(Y)}\RC,
\]
so $d\Psi$ is surjective at $\Gamma_U(\Psi(Y))=\gamma Y$, $\gamma\in\GH$ (Proposition~\ref{prop:BEsection}(3)). By Lemma~\ref{lem:GH}(1), $d\Psi_Y=d\Psi_{\gamma Y}\circ d\gamma_Y$, and $d\gamma_Y\colon T_Y\SC\to T_{\gamma Y}\SC$ is an isomorphism since $\gamma$ is an automorphism of $\SC$; hence $d\Psi_Y$ is surjective.

(2) By Lemma~\ref{lem:GH}(1), the two morphisms $\GH\times\SC\to\RC$, $(\gamma,Y')\mapsto\Psi(\gamma Y')$ and $(\gamma,Y')\mapsto\Psi(Y')$, agree on closed points, and hence agree on $\GH\times\SC$ since it is reduced. Differentiating at $(1,Y)$ along $\GH$ gives $d\Psi\circ c=0$, i.e.\ $\Im c\subset\ker d\Psi$. By (1) and Lemma~\ref{lem:SC},
\[
\dim\ker d\Psi=h^0(\NY(-H))-\dim\RC=(M+1)+(2b-e-2)-(2b-e-2)=M+1=\dim\Im c ,
\]
so $\ker d\Psi=\Im c$.
\end{proof}

Finally, we conclude that the map $a$ is surjective.

\begin{proposition}\label{prop:a-surj}
The map $\bar\delta^{-1}\circ a=d\Psi\colon H^0(\NY(-H))\to H^0(N_C\otimes K_Y|_C)/k\,s_C$ is surjective. Consequently $a$ and $F$ are surjective, and $\ker a=\Im c$.
\end{proposition}

\begin{proof}
By Proposition~\ref{prop:dPsi}, $d\Psi$ is surjective with kernel $\Im c$. By Proposition~\ref{prop:a=dPsi}, $a=\bar\delta\circ d\Psi$ is surjective with kernel $\Im c$. Surjectivity of $F$ follows from diagram \eqref{diag:a} and the sequence \eqref{eq:3.61}, as explained in the proof of Theorem~\ref{thm:alpha}.
\end{proof}

%\begin{proof}[Proof of Theorem~\ref{thm:alpha}]
%By \eqref{eq:alpha-F}, $\alpha(\Yt)=h^0(-H-2K_Y)+\dim\coker F$, and $F$ is surjective by Proposition~\ref{prop:a-surj}. For the last assertion note that $-H-2K_Y=3C_0+(2e+4-b)f$ and $h^0(\OO_Y(3C_0+kf))=\sum_{i=0}^{3}h^0(\OO_{\PP^1}(k-ie))$, which vanishes if and only if $k<0$, i.e.\ $b\ge 2e+5$.
%\end{proof}

%\begin{remark}[Deformation-theoretic meaning]
%The equality $\ker a=\Im c$ obtained in Proposition~\ref{prop:a-surj} says: a first-order deformation $Y_\varepsilon$ of the scroll inside $\PM$ which contains $C_\Dd$ is contained in a first-order deformation $\Yt_\varepsilon$ of the carpet with $\Yt_\varepsilon\cap H_\Dd=\Ct_\Dd$ if and only if $Y_\varepsilon=\gamma_\varepsilon Y$ for a projectivity $\gamma_\varepsilon\in\GH(\Dd)$ fixing $H$ pointwise; and then $\Yt_\varepsilon=\gamma_\varepsilon\Yt$. Equivalently, by Proposition~\ref{prop:BEsection}(2): the hyperplane section of the K3 carpet on a scroll $Y'\supset C$ is the ribbon whose restricted cotangent bundle is $V_{Y'}=E_{Y'}\otimes L_C$ ($\cong E_{Y'}^{\vee}\otimes K_C$), and this ribbon determines $Y'$ up to $\GH$. The scalar $\yt_h$ is the ``horizontal component'' of the carpet class of \cite[Proposition~1.7]{GP97}.
%\end{remark}

\section{Prime K3 surfaces of genus $g \geq 6$ and Fano threefolds of Picard rank one}\label{sec:fano}

In this section we give sharp bounds to $\alpha(X)$, where $X$ is a prime $K3$ surface of genus $g \geq 6$.

\subsection{Set-up}\label{subsec:fano-setup}

For a smooth K3 surface $X\subset\PP^M$ we write $N_X:=N_{X/\PP^M}$ and, as in \S\ref{subsec:setup}, $\alpha(X):=h^0(N_X(-1))-(M+1)$. Recall \cite{Zak91,Lvo92} that if $\alpha(X)\le0$ then $X$ is not extendable: every $V\subset\PP^{M+1}$ having $X$ as a hyperplane section is a cone over $X$.

Fix $g\ge3$. Let $\mathcal H_g$ be the locus, in the Hilbert scheme of $\PP^g$, of smooth K3 surfaces $X\subset\PP^g$ embedded by the complete linear series of a very ample line bundle $H$ with $H^2=2g-2$ and $H$ primitive in $\operatorname{Pic}X$. We call them \emph{prime K3 surfaces of genus $g$}. Since the moduli space of primitively polarized K3 surfaces of genus $g$ is irreducible \cite[Ch.~6]{Huy16}, $\mathcal H_g$ is irreducible. It is smooth of dimension $18+(g+1)^2$, because for $X\in\mathcal H_g$ the Euler sequence gives $h^0(T_{\PP^g}|_X)=(g+1)^2-1$ and $H^1(T_{\PP^g}|_X)=H^2(\OO_X)=k$, the map $H^1(T_X)\to H^2(\OO_X)$, cup product with $c_1(H)$, is surjective, and $H^2(T_X)=0$, so that $h^1(N_X)=0$ and $h^0(N_X)=(g+1)^2-1+19$. A point of $\mathcal H_g$ is called \emph{general} if it lies in a suitable dense open subset.

Let $e\in\{0,1\}$ and $b\ge e+2$, and let $\Yt\subset\PP^g$ be the K3 carpet on $Y=F_e$ embedded by $|H|$, $H=C_0+bf$, so that $g=2b-e+1$ (this is the set-up of \S\ref{subsec:setup}, except that we allow $b=e+2$). By Theorem~\ref{degeneration into $K3$ carpets}, (see also \cite[Theorems~3.5 and 4.4]{GP97}), $[\Yt]$ lies in the closure of $\mathcal H_g$. Since $\Yt$ and the smooth K3 surfaces are local complete intersections, the normal sheaf of a flat family in $\PP^g$ degenerating a general $X\in\mathcal H_g$ to $\Yt$ is locally free and commutes with base change; by upper semicontinuity, for $X$ general in $\mathcal H_g$ and every $k\ge0$,
\begin{equation}\label{eq:semicont}
h^0(N_X(-k))\le h^0(\NYt(-k\Ht)),\qquad\text{in particular}\qquad \alpha(X)\le\alpha(\Yt).
\end{equation}

\subsection{The invariant \texorpdfstring{$\alpha$}{alpha} for \texorpdfstring{$g \geq 6$}{6 <= g <= 12}}

\begin{theorem}\label{thm:alpha-lowgenus}
Let $X\subset\PP^g$ be a general prime K3 surface of genus $g$.
\begin{enumerate}
\item If $g=12$, then $\alpha(X)\le1$.
\item If $g=10$, then $\alpha(X)\le3$.
\item If $g=9$, then $\alpha(X)\le4$.
\item If $g=8$, then $\alpha(X)\le6$.
\item If $g=7$, then $\alpha(X)\le8$.
\item If $g=6$, then $\alpha(X)\le10$.
\item If $g=11$ or $g\ge13$, then $\alpha(X)\le0$; consequently $X$ is not extendable.
\end{enumerate}
\end{theorem}

\begin{proof}
We can degenerate such a K3 surface to a K3 carpet $\Yt\subset\PP^g$, extending $Y\subset\PP^g$, where $Y=F_e$ is embedded by the complete linear series of $C_0+bf$ with $e=0$ or $e=1$ depending on whether $g$ is odd or even respectively, and $g=2b-e+1$. By \eqref{eq:semicont} it is enough to find an upper bound for $\alpha(\Yt)$. In cases (1)--(5) the pair $(e,b)$ is $(1,6)$, $(1,5)$, $(0,4)$, $(1,4)$, $(0,3)$ respectively, and in case (7) it is $(0,b)$ with $b\ge5$ or $(1,b)$ with $b\ge7$. In all these cases $b\ge e+3$, so by Theorem~\ref{thm:alpha}
\[
\alpha(\Yt)=h^0\big(\OO_Y(3C_0+(2e+4-b)f)\big)=\sum_{i=0}^{3}h^0\big(\OO_{\PP^1}(2e+4-b-ie)\big),
\]
which equals $1,3,4,6,8$ in cases (1)--(5) and $0$ in case (7). In case (6), $(e,b)=(1,3)$ and $\alpha(\Yt)=10$ by Proposition~\ref{prop:g6} below. The last assertion in (7) is the theorem of Zak and L'vovsky.
\end{proof}

\begin{proposition}[The case $g=6$]\label{prop:g6}
Let $Y=F_1\subset\PP^6$ be embedded by $|C_0+3f|$ and let $\Yt$ be the K3 carpet on $Y$. Then $\alpha(\Yt)=h^0(\OO_Y(-H-2K_Y))=h^0(\OO_Y(3C_0+3f))=10$.
\end{proposition}

\begin{proof}
Here $b=e+2$, $M=6$, and $\OO_Y(-H-K_Y)=\OO_Y(C_0)$, with $h^0(\OO_Y(C_0))=1$ and $h^1(\OO_Y(C_0))=h^2(\OO_Y(C_0))=0$. We go through the proof of Theorem~\ref{thm:alpha}, where the assumption $b\ge e+3$ was used only to have $H^0(-H-K_Y)=0$ in \eqref{eq:3.41} and \eqref{eq:3.52}; everything else, in particular \eqref{eq:3.47}, \eqref{eq:3.48}, \eqref{eq:3.56}--\eqref{eq:3.59} and \eqref{eq:3.61}, holds for $b\ge e+1$. From \eqref{eq:3.41} and \eqref{eq:3.48} we get an exact sequence
\begin{multline*}
0\to H^0(\NYt(-H+K_Y))\to H^0(\OO_Y(C_0))=k\xrightarrow{\ \delta_1\ }H^1(\NY(-H+K_Y))=k^3\\
\xrightarrow{\ j\ }H^1(\NYt(-H+K_Y))\to0
\end{multline*}
and $H^2(\NYt(-H+K_Y))=0$; thus $h^0(\NYt(-H+K_Y))=1-\operatorname{rank}\delta_1$, $h^1(\NYt(-H+K_Y))=3-\operatorname{rank}\delta_1$, and $j$ is surjective. From \eqref{eq:3.52} and \eqref{eq:3.58}, $h^0(\Hc(-H))=h^0(\NY(-H))-1=(M+1)+h^1(T_Y(-H))-1=7+3-1=9$, and from \eqref{eq:3.61}, $h^0(\NYt(-H))=9+h^0(\OO_Y(3C_0+3f))=19$ and $H^1(\NYt(-H))=H^1(\OO_Y(3C_0+3f))=0$. Sequence \eqref{eq:3.40} therefore gives
\[
h^0(\NYt(-\Ht))=h^0(\NYt(-H+K_Y))+h^0(\NYt(-H))-\operatorname{rank}F=20-\operatorname{rank}\delta_1-\operatorname{rank}F,
\]
where $F\colon H^0(\NYt(-H))\to H^1(\NYt(-H+K_Y))$ is the connecting map, and it remains to show that $F$ is surjective. Since then, $\operatorname{rank}F=3-\operatorname{rank}\delta_1$ and $\alpha(\Yt)=17-7=10$.

The proof of Proposition~\ref{prop:a=dPsi} shows, without using $b\ge e+3$, that $F(\iota(x))=j(\delta(t))$ for every $x\in H^0(\NY(-H))$, where $\iota\colon H^0(\NY(-H))\to H^0(\NYt(-H))$ is induced by restriction of homomorphisms from $I_Y$ to $I_{\Yt}$, $d\Psi(x)=[t]$, and $\delta$ is the map of Lemma~\ref{lem:delta}. Now $\delta$ is surjective (Lemma~\ref{lem:delta}), $d\Psi$ is surjective (Proposition~\ref{prop:dPsi}; its proof, and those of the results of \S\S\ref{subsec:geometric}--\ref{subsec:sections} on which it relies, only use $b\ge e+1$), and $j$ is surjective; hence the image of $F$ is all of $H^1(\NYt(-H+K_Y))$.
\end{proof}

\subsection{Sections of \texorpdfstring{$N(-2)$}{N(-2)}}

In this subsection we give upper bounds to the dimension of the space $H^0(N_{\widetilde Y}(-2))$ where $\widetilde Y \subset \mathbb{P}^g$ is a $K3$ carpet on a rational normal scroll $Y \subset \mathbb{P}^g$. This consequently gives upper bounds to the dimension of the space $H^0(N_{X}(-2))$, where $X \subset \mathbb{P}^g$ is a prime $K3$ surface of genus $g$. 

\begin{proposition}\label{prop:N2}
Let $Y=F_e\subset\PM$ be embedded by $|H|$, $H=C_0+bf$ with $b\ge e+2$, and let $\Yt$ be the K3 carpet on $Y$. Then
\[
h^0(\NYt(-2\Ht))\le h^0(\OO_Y(-2H-2K_Y))=h^0\big(\OO_Y(2C_0+(2e+4-2b)f)\big).
\]
In particular $h^0(\NYt(-2\Ht))=0$ if $b\ge e+3$, and $h^0(\NYt(-2\Ht))\le1$ if $(e,b)=(1,3)$.
\end{proposition}

\begin{proof}
Tensoring \eqref{eq:3.1} by $\OO_{\Yt}(-2\Ht)$ gives $h^0(\NYt(-2\Ht))\le h^0(\NYt(-2H+K_Y))+h^0(\NYt(-2H))$. We use throughout that $H^0(\OO_Y(-H+K_Y))=0$, $H^0(\OO_Y(-2H-K_Y))=0$ (as $-2H-K_Y=(e+2-2b)f$ and $2b>e+2$), and, by Kodaira vanishing, $H^1(\OO_Y(-2H))=H^1(\OO_Y(-2H+K_Y))=0$ and $H^i(\OO_Y(-H))=0$ for all $i$ (see the derivation of \eqref{eq:3.54}).

\emph{Step 1: $H^0(\NYt(-2H+K_Y))=0$.} By \eqref{eq:3.2} tensored with $-2H+K_Y$, $H^0(\NYt(-2H+K_Y))=H^0(\Hc(-2H+K_Y))$, and by \eqref{eq:3.3} tensored with $-2H+K_Y$ the latter is a quotient of $H^0(\NY(-2H+K_Y))$. By \eqref{eq:3.11} tensored with $-2H+K_Y$, $H^0(\TP(-2H+K_Y))=0$, so by \eqref{eq:3.10} $H^0(\NY(-2H+K_Y))$ injects into $H^1(T_Y(-2H+K_Y))$. Finally $H^1(T_Y(-2H+K_Y))=0$: in
\[
0\to T_{Y/\PP^1}(-2H+K_Y)\to T_Y(-2H+K_Y)\to p^*T_{\PP^1}(-2H+K_Y)\to0
\]
we have $T_{Y/\PP^1}(-2H+K_Y)=\OO_Y(-2C_0-(2b+2)f)$ and $p^*T_{\PP^1}(-2H+K_Y)=\OO_Y(-4C_0-(2b+e)f)$; both have vanishing $p_*$, while by relative duality $R^1p_*\OO_Y(-2C_0)=\OO_{\PP^1}(e)$ and $R^1p_*\OO_Y(-4C_0)=\OO_{\PP^1}(e)\oplus\OO_{\PP^1}(2e)\oplus\OO_{\PP^1}(3e)$. Hence
\begin{gather*}
H^1(T_{Y/\PP^1}(-2H+K_Y))=H^0(\OO_{\PP^1}(e-2b-2))=0,\\
H^1(p^*T_{\PP^1}(-2H+K_Y))=H^0\big(\OO_{\PP^1}(-2b)\oplus\OO_{\PP^1}(e-2b)\oplus\OO_{\PP^1}(2e-2b)\big)=0 .
\end{gather*}

\emph{Step 2: $H^0(\Hc(-2H))=0$.} Granting this, \eqref{eq:3.2} tensored with $-2H$ gives $h^0(\NYt(-2H))\le h^0(\OO_Y(-2H-2K_Y))$, which together with Step~1 proves the inequality. By \eqref{eq:3.3} tensored with $-2H$ we have an exact sequence
\[
0\to H^0(\NY(-2H))\to H^0(\Hc(-2H))\to H^1(\OO_Y(-2H-K_Y))\xrightarrow{\ s^{\vee}\ }H^1(\NY(-2H)),
\]
so it suffices to show that (a) $H^0(\NY(-2H))=0$ and that (b) $s^{\vee}$ is injective.

(a) Twisting the sequence \eqref{eq:3.10} by $-2H$ and taking cohomology, we have
\[
0\to H^0(T_Y(-2H))\to H^0(\TP(-2H))\to H^0(\NY(-2H))\to H^1(T_Y(-2H))\xrightarrow{\ s\ }H^1(\TP(-2H)) .
\]
Turning to the Euler sequence \eqref{eq:3.11}, twisting by $-2H$ and taking cohomology, we get
\[
0\to H^0(\OO_Y(-H)^{\oplus M+1})=(0)\to H^0(\TP(-2H))\to H^1(\OO_Y(-2H))=(0) .
\]
Hence $H^0(\TP(-2H))=0$. So it is enough to show that the map $H^1(T_Y(-2H))\xrightarrow{\ s\ }H^1(\TP(-2H))$ is injective.

First consider the two exact sequences, the first one is the relative Euler sequence on $Y=\PP(\OO_{\PP^1}\oplus\OO_{\PP^1}(-e))$ and the second one is the Euler sequence of $\PM$ pulled back to $Y$.
\[
\begin{tikzcd}
  0 \arrow[r] & \OO_Y \arrow[r] & p^*(\OO_{\PP^1} \oplus \OO_{\PP^1}(-e))^* \otimes \OO_{Y}(C_0) \arrow[r, "x"] \arrow[dr, "y \circ x"] & T_{Y/\PP^1} \arrow[r] \arrow[d, "y"] & 0 \\
   0 \arrow[r] & \OO_Y \arrow[r] & \OO_Y(H)^{\oplus M+1} \arrow[r] & \TP \arrow[r] & 0
\end{tikzcd}
\]
The map $y$ is the injective map induced by the composition of two injective maps given by $0\to T_{Y/\PP^1}\to T_Y$ and $0\to T_Y\to\TP$. We show that the map $y$ has a unique lift to a map
\[
z\colon p^*(\OO_{\PP^1}\oplus\OO_{\PP^1}(-e))^*\otimes\OO_Y(C_0)\to\OO_Y(H)^{\oplus M+1} .
\]
Apply the functor $\Hom(p^*(\OO_{\PP^1}\oplus\OO_{\PP^1}(-e))^*\otimes\OO_Y(C_0),-)$ to the bottom exact sequence we have
\begin{multline*}
0\to\Hom(p^*(\OO_{\PP^1}\oplus\OO_{\PP^1}(-e))^*\otimes\OO_Y(C_0),\OO_Y)\to\Hom(p^*(\OO_{\PP^1}\oplus\OO_{\PP^1}(-e))^*\otimes\OO_Y(C_0),\OO_Y(H)^{\oplus M+1})\\
\to\Hom(p^*(\OO_{\PP^1}\oplus\OO_{\PP^1}(-e))^*\otimes\OO_Y(C_0),\TP)\to\Ext^1(p^*(\OO_{\PP^1}\oplus\OO_{\PP^1}(-e))^*\otimes\OO_Y(C_0),\OO_Y)
\end{multline*}
Now, $\Ext^1(p^*(\OO_{\PP^1}\oplus\OO_{\PP^1}(-e))^*\otimes\OO_Y(C_0),\OO_Y)=\Ext^1(\OO_Y(C_0)\oplus\OO_Y(C_0+ef),\OO_Y)=H^1(\OO_Y(-C_0)\oplus\OO_Y(-C_0-ef))=0$. Further, $\Hom(p^*(\OO_{\PP^1}\oplus\OO_{\PP^1}(-e))^*\otimes\OO_Y(C_0),\OO_Y)=\Hom(\OO_Y(C_0)\oplus\OO_Y(C_0+ef),\OO_Y)=H^0(\OO_Y(-C_0)\oplus\OO_Y(-C_0-ef))=0$. Hence the map $y$ has a unique lift to the map $z$, which induces the following commutative diagram involving the two previous exact sequences
\begin{equation}\label{diag:relEuler}
\begin{tikzcd}
  0 \arrow[r] & \OO_Y \arrow[r] \arrow[d, "w", "\cong"'] & p^*(\OO_{\PP^1} \oplus \OO_{\PP^1}(-e))^* \otimes \OO_Y(C_0) \arrow[r, "x"] \arrow[d, "z"] & T_{Y/\PP^1} \arrow[r] \arrow[d, "y"] & 0 \\
   0 \arrow[r] & \OO_Y \arrow[r] & \OO_Y(H)^{\oplus M+1} \arrow[r] & \TP \arrow[r] & 0
\end{tikzcd}
\end{equation}
To see that the map $w$ is an isomorphism, we only need to show that $w$ is non-zero. To show this restrict the above diagram to a fiber of the projection map $p\colon Y\to\PP^1$. Since $\OO_Y(H)=\OO_Y(C_0+bf)$, we have that under the embedding, the fiber $\PP^1$ is embedded by the complete linear series of $\OO_{\PP^1}(1)$, i.e.\ by a linear embedding $i\colon\PP^1\hookrightarrow\PM$. Hence after restricting the above diagram to the fiber, we get
\[
\begin{tikzcd}
  0 \arrow[r] & \OO_{\PP^1} \arrow[r, "\sigma_{\PP^1}"] \arrow[d, "\Bar{w}", "\cong"'] & \OO_{\PP^1}(1)^{\oplus 2} \arrow[r, "\rho_{\PP^1}"] \arrow[d, "\widetilde{\alpha}"] & T_{\PP^1} \arrow[r] \arrow[d, "\beta = di"] & 0 \\
   0 \arrow[r] & \OO_{\PP^1} \arrow[r, "i^*\sigma_{\PM}"] & \OO_{\PP^1}(1)^{\oplus M+1} \arrow[r, "i^*\rho_{\PM}"] & T_{\PP^1} \oplus \OO_{\PP^1}(1)^{\oplus M-1} = T_{\PM}|_{\PP^1} \arrow[r] & 0
\end{tikzcd}
\]
where $\sigma_{\PP^1}$ and $\rho_{\PP^1}$ are the maps of the Euler exact sequence on $\PP^1$, while $i^*\sigma_{\PM}$ and $i^*\rho_{\PM}$ are the maps of the Euler exact sequence on $\PM$, pulled back to the line $\PP^1$, and $\beta=di$ is the differential map of the embedding $i$. As before, since $\Hom(\OO_{\PP^1}(1)^{\oplus 2},\OO_{\PP^1})=\Ext^1(\OO_{\PP^1}(1)^{\oplus 2},\OO_{\PP^1})=0$, the map $\widetilde{\alpha}$ is the unique map lifting $di\circ\rho_{\PP^1}$. Choose coordinates $[u:v]$ on $\PP^1$ and $[X_0:\dots:X_M]$ on $\PM$ such that $i([u:v])=[u:v:0:\dots:0]$. Then $di$ is the inclusion of the first summand $T_{\PP^1}\to T_{\PP^1}\oplus\OO_{\PP^1}(1)^{\oplus M-1}$, and the map
\[
\widetilde\alpha\colon\OO_{\PP^1}(1)^{\oplus 2}\longrightarrow\OO_{\PP^1}(1)^{\oplus M+1},\qquad \widetilde\alpha(a,b)=(a,b,0,\dots,0),
\]
makes the diagram commutative, since both $\rho_{\PP^1}(a,b)$ and $i^*\rho_{\PM}(a,b,0,\dots,0)$ are the vector field $a\,\partial/\partial X_0+b\,\partial/\partial X_1$ along the line. Hence this is $\widetilde\alpha$, and $\widetilde{\alpha}\circ\sigma_{\PP^1}(1)=\widetilde{\alpha}(u,v)=(u,v,0,\dots,0)=i^*\sigma_{\PM}(1)$. Therefore $\Bar{w}(1)=1$. Hence $w$ is an isomorphism.

Now twisting the commutative diagram \eqref{diag:relEuler} by $\OO_Y(-2H)$, we have the following
\begin{equation}\label{diag:relEuler2}
\begin{tikzcd}
  0 \arrow[r] & \OO_Y(-2H) \arrow[r] \arrow[d, "w(-2H)", "\cong"'] & p^*(\OO_{\PP^1} \oplus \OO_{\PP^1}(-e))^* \otimes \OO_Y(C_0-2H) \arrow[r, "x"] \arrow[d, "z"] & T_{Y/\PP^1}(-2H) \arrow[r] \arrow[d, "y"] & 0 \\
   0 \arrow[r] & \OO_Y(-2H) \arrow[r] & \OO_Y(-H)^{\oplus M+1} \arrow[r] & \TP(-2H) \arrow[r] & 0
\end{tikzcd}
\end{equation}
Taking cohomology of the bottom row of diagram \eqref{diag:relEuler2}, since $H^1(\OO_Y(-H))=H^2(\OO_Y(-H))=0$, we get
\[
0\to H^1(\TP(-2H))\to H^2(\OO_Y(-2H))\to0 .
\]
Similarly tensoring the exact sequence of horizontal and vertical tangent bundles of $Y$ by $-2H$, we have
\[
0\to T_{Y/\PP^1}(-2H)\to T_Y(-2H)\to p^*(T_{\PP^1})(-2H)\to0
\]
and taking the cohomology, we have that
\[
H^0(p^*(T_{\PP^1})(-2H))\to H^1(T_{Y/\PP^1}(-2H))\to H^1(T_Y(-2H))\to H^1(p^*(T_{\PP^1})(-2H))\to0 .
\]
Now note that $p^*(T_{\PP^1})(-2H)=-2C_0+(2-2b)f$. We have that $H^0(-2C_0+(2-2b)f)=0$ and $H^1(-2C_0+(2-2b)f)=H^1((2b-e-4)f)^*=0$ if $b\ge e+2$. Hence we have
\[
0\to H^1(T_{Y/\PP^1}(-2H))\to H^1(T_Y(-2H))\to0 .
\]
On the other hand, $p^*(\OO_{\PP^1}\oplus\OO_{\PP^1}(-e))^*\otimes\OO_Y(C_0-2H)=\OO_Y(-C_0-2bf)\oplus\OO_Y(-C_0+(e-2b)f)$ has no cohomology, since $p_*\OO_Y(-C_0)=R^1p_*\OO_Y(-C_0)=0$; so taking the cohomology of the top row of the diagram \eqref{diag:relEuler2}, we get
\[
0\to H^1(T_{Y/\PP^1}(-2H))\to H^2(\OO_Y(-2H))\to0 .
\]
Hence taking the cohomology of the diagram \eqref{diag:relEuler2}, we have a commutative diagram
\[
\begin{tikzcd}
    H^1(T_{Y/\PP^1}(-2H)) \arrow[d, "\cong"] \arrow[r, "\cong"] & H^2(\OO_Y(-2H)) \arrow[dd, "\cong"] \\
     H^1(T_{Y}(-2H)) \arrow[d, "s"] & \\
     H^1(\TP(-2H)) \arrow[r, "\cong"] & H^2(\OO_Y(-2H))
\end{tikzcd}
\]
This shows that the map $s$ is injective and is in fact an isomorphism. Hence $H^0(\NY(-2H))=0$.
\smallskip

(b) Consider the exact sequence given by twisting \eqref{eq:3.3} by $-2H$,
\begin{equation}\label{eq:b16}
0\to\OO_Y(-2H-K_Y)\xrightarrow{\ s^{\vee}\ }\NY(-2H)\to\Hc(-2H)\to0 .
\end{equation}
Let $p\colon Y\to\PP^1$ denote the smooth fibration. Now note that $-2H-K_Y=(e+2-2b)f$ and hence $R^1p_*(-2H-K_Y)=0$. Hence we have the following exact sequence by pushing forward the exact sequence \eqref{eq:b16}
\begin{equation}\label{eq:b17}
0\to p_*(-2H-K_Y)\xrightarrow{\ \beta\ }p_*(\NY(-2H))\to p_*(\Hc(-2H))\to0 .
\end{equation}
It is enough to show that the exact sequence \eqref{eq:b17} is split: then $H^1(p_*(-2H-K_Y))\to H^1(p_*(\NY(-2H)))$ is injective, and since $H^1(\OO_Y(-2H-K_Y))=H^1(p_*(-2H-K_Y))$ and $H^1(p_*(\NY(-2H)))\to H^1(\NY(-2H))$ is injective by the Leray spectral sequence, the map $s^{\vee}\colon H^1(\OO_Y(-2H-K_Y))\to H^1(\NY(-2H))$ is injective. To see this we first show that there is a surjection from $p_*(\NY(-2H))$ to $p_*(-2H-K_Y)$. To see this consider the exact sequence \eqref{eq:3.10} twisted by $-2H$ to get
\begin{equation}\label{eq:b18}
0\to T_Y(-2H)\to\TP(-2H)\to\NY(-2H)\to0 .
\end{equation}
We show that
\[
R^1p_*(\TP(-2H))=0,\qquad p_*(T_Y(-2H))=\OO_{\PP^1}(e-2b)\quad\text{and}\quad R^1p_*(T_Y(-2H))=\OO_{\PP^1}(e+2-2b) .
\]
To see the first vanishing once again note that from the Euler exact sequence \eqref{eq:3.11} twisted by $-2H$, we have
\[
0\to\OO_Y(-2H)\to\OO_Y(-H)^{\oplus M+1}\to\TP(-2H)\to0 .
\]
Since $R^1p_*(\OO_Y(-H))=R^2p_*(\OO_Y(-2H))=0$, we have that $R^1p_*(\TP(-2H))=0$. Now we turn to $R^ip_*(T_Y(-2H))$. Consider the exact sequence of horizontal and vertical tangent bundles on $Y$ twisted by $-2H$
\begin{equation}\label{eq:b19}
0\to T_{Y/\PP^1}(-2H)\to T_Y(-2H)\to p^*T_{\PP^1}(-2H)\to0 .
\end{equation}
Using the relative Euler exact sequence we have that $T_{Y/\PP^1}=2C_0+ef$, so that $T_{Y/\PP^1}(-2H)=(e-2b)f$. Therefore $p_*(T_{Y/\PP^1}(-2H))=\OO_{\PP^1}(e-2b)$ and $R^ip_*(T_{Y/\PP^1}(-2H))=R^ip_*((e-2b)f)=0$ for $i=1,2$. Now $p^*T_{\PP^1}(-2H)=-2C_0+(2-2b)f$. Therefore, $p_*p^*T_{\PP^1}(-2H)=p_*(-2C_0+(2-2b)f)=\OO_{\PP^1}(2-2b)\otimes p_*(-2C_0)=0$ and
\begin{align*}
R^1p_*(p^*T_{\PP^1}(-2H))&=R^1p_*(-2C_0+(2-2b)f)\\
&=\OO_{\PP^1}(2-2b)\otimes R^1p_*(-2C_0)\\
&=\OO_{\PP^1}(2-2b)\otimes(p_*(\OO_Y))^*\otimes(\textstyle\bigwedge^2(\OO_{\PP^1}\oplus\OO_{\PP^1}(-e)))^*\\
&=\OO_{\PP^1}(e+2-2b) .
\end{align*}
Pushing forward the sequence \eqref{eq:b19}, we have that
\[
0\to p_*(T_{Y/\PP^1}(-2H))=\OO_{\PP^1}(e-2b)\to p_*T_Y(-2H)\to0
\]
and
\[
0\to R^1p_*T_Y(-2H)\to\OO_{\PP^1}(e+2-2b)\to0 .
\]
Now pushing forward the exact sequence \eqref{eq:b18}, we have that
\begin{multline}\label{eq:b23}
0\to\OO_{\PP^1}(e-2b)=p_*T_Y(-2H)=p_*(T_{Y/\PP^1}(-2H))\to p_*(\TP(-2H))\to p_*(\NY(-2H))\\
\xrightarrow{\ \lambda\ }\OO_{\PP^1}(e+2-2b)=p_*(-2H-K_Y)\to0 .
\end{multline}
We now end the proof by showing that for the K3 carpet $\Yt$, which is defined by the injection $s^{\vee}\colon\OO_Y(-K_Y)\to\NY$ of \eqref{eq:s}, the map $\lambda$ is, up to a nonzero scalar, a section to the injective map $\beta$ induced by $s^{\vee}$ after twisting by $-2H$ and pushing forward, which sits in the following diagram
\begin{equation}\label{diag:section}
\begin{tikzcd}
  &&& 0 \arrow[d] && \\
  &&& p_*(-2H-K_Y) \arrow[d, "\beta"] && \\
  0 \arrow[r] & \OO_{\PP^1}(e-2b) \arrow[r] &  p_*(\TP(-2H)) \arrow[r, "\delta"] & p_*(\NY(-2H)) \arrow[r, "\lambda"] & p_*(-2H-K_Y) \arrow[r] & 0
\end{tikzcd}
\end{equation}
Since $p_*(-2H-K_Y)=\OO_{\PP^1}(e+2-2b)$ is a line bundle, it is enough to show that $\lambda\circ\beta$ is non-zero. This happens if and only if $\operatorname{Im}(\beta)$ is not contained in $\operatorname{Ker}(\lambda)$. So we need to show that $\beta$ does not factor as
\[
\begin{tikzcd}
   & & p_*(-2H-K_Y) \arrow[d, "\beta"] \arrow[dl] \\
   0 \arrow[r] & \operatorname{Ker}(\lambda) \arrow[r] & p_*(\NY(-2H))
\end{tikzcd}
\]
i.e.\ that $\beta$ is not in the image of
\begin{equation}\label{eq:b32}
\Hom(p_*(-2H-K_Y),\operatorname{Ker}(\lambda))\to\Hom(p_*(-2H-K_Y),p_*(\NY(-2H))) .
\end{equation}
To see this note that (by sequence \eqref{eq:b23} and diagram \eqref{diag:section}), we have an exact sequence
\[
0\to\OO_{\PP^1}(e-2b)=p_*T_{Y/\PP^1}(-2H)\to p_*(\TP(-2H))\to\operatorname{Ker}(\lambda)\to0 .
\]
Applying the functor $\Hom(p_*(-2H-K_Y),-)$ we get
\begin{multline}\label{eq:b30}
0\to\Hom(p_*(-2H-K_Y),p_*T_{Y/\PP^1}(-2H))\to\Hom(p_*(-2H-K_Y),p_*(\TP(-2H)))\\
\to\Hom(p_*(-2H-K_Y),\operatorname{Ker}(\lambda))\to\Ext^1(p_*(-2H-K_Y),p_*T_{Y/\PP^1}(-2H))\\
\xrightarrow{\ \eta\ }\Ext^1(p_*(-2H-K_Y),p_*(\TP(-2H))) .
\end{multline}
We first show that the map $\eta$ is an injective map and is in fact an isomorphism. Pushing forward the bottom row of diagram \eqref{diag:relEuler2}, since $p_*(\OO_Y(-H))=R^1p_*(\OO_Y(-H))=0$, we get
\begin{equation}\label{eq:b26}
0\to p_*(\TP(-2H))\to R^1p_*(-2H)\to0 ,
\end{equation}
and pushing forward the top row of diagram \eqref{diag:relEuler2}, since $p^*(\OO_{\PP^1}\oplus\OO_{\PP^1}(-e))^*\otimes\OO_Y(C_0-2H)$ has no direct images, we get
\begin{equation}\label{eq:b29}
0\to p_*T_{Y/\PP^1}(-2H)\to R^1p_*(-2H)\to0 .
\end{equation}
Because of the commutative diagram \eqref{diag:relEuler2}, \eqref{eq:b26} and \eqref{eq:b29} sit in a commutative diagram
\[
\begin{tikzcd}
0 \arrow[r] & \Ext^1(p_*(-2H-K_Y), p_*T_{Y/\PP^1}(-2H)) \arrow[d, "\eta"] \arrow[r] & \Ext^1(p_*(-2H-K_Y), R^1p_*(-2H)) \arrow[r] \arrow[d, "\cong"] & 0 \\
0 \arrow[r] & \Ext^1(p_*(-2H-K_Y), p_*(\TP(-2H))) \arrow[r] & \Ext^1(p_*(-2H-K_Y), R^1p_*(-2H)) \arrow[r] & 0
\end{tikzcd}
\]
The right hand vertical map above is an isomorphism because it is induced by $w(-2H)$ in diagram \eqref{diag:relEuler2} which is an isomorphism. Hence the map $\eta$ is an isomorphism. Going back to sequence \eqref{eq:b30}, we have
\begin{multline}\label{eq:b31}
0\to\Hom(p_*(-2H-K_Y),p_*T_{Y/\PP^1}(-2H))\to\Hom(p_*(-2H-K_Y),p_*(\TP(-2H)))\\
\to\Hom(p_*(-2H-K_Y),\operatorname{Ker}(\lambda))\to0 .
\end{multline}
Therefore, the image of the map \eqref{eq:b32} is the image of
\begin{equation}\label{eq:b33}
\Hom(p_*(-2H-K_Y),p_*(\TP(-2H)))\to\Hom(p_*(-2H-K_Y),p_*(\NY(-2H))) .
\end{equation}
Since $p^*p_*(-2H-K_Y)=-2H-K_Y$, by the projection formula, this is the map
\begin{equation}\label{eq:b34}
H^0(\TP\otimes K_Y)\to H^0(\NY\otimes K_Y),
\end{equation}
and $\beta$ corresponds to $s\in H^0(\NY\otimes K_Y)$. Now by taking the cohomology of the exact sequence
\[
0\to T_Y\otimes K_Y\to\TP\otimes K_Y\to\NY\otimes K_Y\to0
\]
we have
\[
0\to H^0(T_Y\otimes K_Y)\to H^0(\TP\otimes K_Y)\to H^0(\NY\otimes K_Y)\xrightarrow{\ \partial\ }H^1(T_Y\otimes K_Y) .
\]
Now $s$ is not in the image of \eqref{eq:b34}, because $\partial(s)\neq0$ by Lemma~\ref{lem:pushforward}. Hence $\lambda\circ\beta\neq0$, the sequence \eqref{eq:b17} is split, and $s^{\vee}\colon H^1(\OO_Y(-2H-K_Y))\to H^1(\NY(-2H))$ is injective. Therefore $H^0(\Hc(-2H))=0$.
\end{proof}

\subsection{Classification of Fano threefolds of Picard rank one}

\begin{definition}\label{def:Vg}
For $g\ge3$ let $\mathcal V_g$ denote the locus, in the Hilbert scheme of $\PP^{g+1}$, of smooth Fano threefolds $V\subset\PP^{g+1}$ of Picard rank one, index one and genus $g$, i.e.\ $\operatorname{Pic}V=\mathbb Z(-K_V)$ and $(-K_V)^3=2g-2$, embedded by the complete linear series $|-K_V|$; it is open in the Hilbert scheme, and $\overline{\mathcal V_g}$ denotes its closure. For $X\subset\PP^g$ we denote by $C(X)\subset\PP^{g+1}$ the cone over $X$.
\end{definition}

\begin{corollary}\label{cor:cone}
Let $X\subset\PP^g$ be a general prime K3 surface of genus $g\ge6$. Then the tangent space to the Hilbert scheme of $\PP^{g+1}$ at $[C(X)]$ has dimension
\[
\dim T_{[C(X)]}=h^0(N_X)+h^0(N_X(-1))+h^0(N_X(-2))=18+(g+1)^2+(g+1)+\alpha(X)+h^0(N_X(-2)),
\]
where $h^0(N_X(-2))=0$ if $g\ge7$ and $h^0(N_X(-2))\le1$ if $g=6$.
\end{corollary}

\begin{proof}
Since $X$ is projectively normal and $H^1(\OO_X(k))=0$ for all $k$, one has $H^0(N_{C(X)/\PP^{g+1}})=\bigoplus_{k\ge0}H^0(N_X(-k))$ (see \cite[\S2]{CLM93}). We have $h^0(N_X)=18+(g+1)^2$ (\S\ref{subsec:fano-setup}) and $h^0(N_X(-1))=(g+1)+\alpha(X)$ by definition. Since $\operatorname{Pic}X=\mathbb ZH$ and $g\ge5$, the homogeneous ideal of $X$ is generated by quadrics \cite{SD74}; hence $N_X=\sHom(I_X/I_X^2,\OO_X)$ is a subsheaf of a direct sum of copies of $\OO_X(2)$, and $H^0(N_X(-k))=0$ for $k\ge3$. Finally $h^0(N_X(-2))\le h^0(\NYt(-2\Ht))$ by \eqref{eq:semicont}, and the right-hand side is $0$ for $g\ge7$ (then $b\ge e+3$) and at most $1$ for $g=6$ ($(e,b)=(1,3)$) by Proposition~\ref{prop:N2}.
\end{proof}

\begin{lemma}\label{lem:hypsection}
Let $W$ be an irreducible component of $\overline{\mathcal V_g}$ and $V$ a general element of $W$. Then a general hyperplane section of $V$ corresponds to a general K3 surface in $\mathcal H_g$.
\end{lemma}

\begin{proof}
The proof follows from standard deformation theory arguments (see for example \cite[Lemma~3.7]{CLM98}).
\end{proof}

\begin{theorem}\label{thm:fano}
Let $g\ge6$.
\begin{enumerate}
\item $\mathcal V_g=\emptyset$ if $g=11$ or $g\ge13$.
\item For $6\le g\le10$ and $g=12$, $\mathcal V_g$ is irreducible, and the families of Fano and Iskovskikh \cite{Isk77,Isk78} (see \cite{Muk88} and \cite[Table~3.1]{CLM98}) form an open dense subset of $\mathcal V_g$.
\end{enumerate}
\end{theorem}

\begin{proof}
\emph{Proof of (1).} If $\mathcal V_g\neq\emptyset$, let $V$ be a general element of an irreducible component of $\overline{\mathcal V_g}$. By Lemma~\ref{lem:hypsection} a general hyperplane section of $V$ is a general prime K3 surface $X$ of genus $g$, and $V$ is not a cone; this contradicts Theorem~\ref{thm:alpha-lowgenus}(7).

\emph{Proof of (2).} First one notes that $\mathcal H_g$ is irreducible. Now the list of examples of Fano--Iskovskikh--Mori--Mukai (see for example \cite[Table~3.1]{CLM98}) shows that in each of the cases mentioned in part (2), there exist K3 surfaces in $\mathcal H_g$ which are extendable to smooth Fano threefolds. Hence by Lemma~\ref{lem:hypsection}, in each of the cases, $\mathcal H_g$ is dominated by $\mathcal V_g$. We end the proof by showing that $\overline{\mathcal V_g}$ has a unique irreducible component and that the examples of Fano--Iskovskikh--Mori--Mukai form a dense subset of $\mathcal V_g$. To see this, we observe first, as in \cite{CLM93} and \cite{CLM98}, that a projectively Cohen--Macaulay scheme degenerates along a flat family to the cone over its hyperplane section. Hence, once again using Lemma~\ref{lem:hypsection}, it is enough to show that for a general K3 surface $X\subset\PP^g$ in $\mathcal H_g$, the cone $C(X)\subset\PP^{g+1}$, which lies on every irreducible component of $\overline{\mathcal V_g}$, is a smooth point of the Hilbert scheme. For the stated values of $g$, Corollary~\ref{cor:cone} along with Theorem~\ref{thm:alpha-lowgenus} gives an upper bound to $\dim T_{[C(X)]}$, which is listed as follows:
\begin{center}
\begin{tabular}{c|cccccc}
$g$ & $6$ & $7$ & $8$ & $9$ & $10$ & $12$\\ \hline
$\dim T_{[C(X)]}\le$ & $85$ & $98$ & $114$ & $132$ & $153$ & $201$
\end{tabular}
\end{center}
Now for a smooth Fano threefold $V\in\mathcal V_g$, $V$ represents a smooth point of the Hilbert scheme, since embedded deformations of an anticanonically embedded Fano threefold are unobstructed by the vanishing of $h^1(N_{V/\PP^{g+1}})$, which follows from the Kodaira--Nakano vanishing theorem. Therefore the dimension of $\overline{\mathcal V_g}$ at $V$ is given by $\dim T_{[V]}=h^0(N_{V/\PP^{g+1}})$. Now using the descriptions of the examples of Fano--Iskovskikh--Mori--Mukai, one can calculate the dimension $h^0(N_{V/\PP^{g+1}})$ (see for example \cite[Table~2, p.~666]{CLM93} and \cite[Table~3.1]{CLM98}) to check that in each case $\dim T_{[V]}$ is exactly the upper bound just obtained. Since $C(X)$ lies on the irreducible component of $\overline{\mathcal V_g}$ containing $V$, whose dimension is $h^0(N_{V/\PP^{g+1}})$, $C(X)$ is also a smooth point of the Hilbert scheme. Hence $\overline{\mathcal V_g}$ has a unique irreducible component. Finally, the examples of Fano--Iskovskikh--Mori--Mukai form an irreducible family (in each case they are the linear sections of a fixed homogeneous variety, see \cite{Muk88}) of dimension $h^0(N_{V/\PP^{g+1}})=\dim\mathcal V_g$, hence a dense subset of $\mathcal V_g$.
\end{proof}

\begin{corollary}\label{cor:equality}
Let $X\subset\PP^g$ be a general prime K3 surface of genus $g$. Then the upper bounds of Theorem~\ref{thm:alpha-lowgenus} and Corollary~\ref{cor:cone} are equalities, namely
\begin{enumerate}
    \item If $g=12$, then $\alpha(X)=1$.
    \item If $g=10$, then $\alpha(X)=3$.
    \item If $g=9$, then $\alpha(X)=4$.
    \item If $g=8$, then $\alpha(X)=6$.
    \item If $g=7$, then $\alpha(X)=8$.
    \item If $g=6$, then $\alpha(X)=10$ and $h^0(N_X(-2))=1$.
\end{enumerate}
\end{corollary}

\begin{proof}
By Corollary~\ref{cor:cone}, $\dim T_{[C(X)]}=18+(g+1)^2+(g+1)+\alpha(X)+h^0(N_X(-2))$, and  $\alpha(X)\le\alpha(\Yt)$ and $h^0(N_X(-2))\le h^0(\NYt(-2\Ht))$ by \eqref{eq:semicont}. Hence $\dim T_{[C(X)]} \leq 18+(g+1)^2+(g+1)+\alpha(\widetilde Y)+h^0(N_{\widetilde Y}(-2\widetilde H)) = h^0(N_{V/\PP^{g+1}})$ by the proof of part $(2)$ of Theorem~\ref{thm:fano}. But a smooth $V$ is a smooth point of the Hilbert scheme and hence $\alpha(X) = \alpha(\Yt)$ and $h^0(N_{X}(-2)) = h^0(N_{\widetilde Y}(-2\widetilde H))$. Now the conclusion follows from Theorem~\ref{thm:alpha}, Proposition~\ref{prop:g6} and Proposition~\ref{prop:N2}.
\end{proof}

\appendix
\section{Extensions of ribbons to \texorpdfstring{$K$}{K}-trivial carpets}\label{sec:extensions}
\begin{center}\textsc{by Anand Deopurkar and Jayan Mukherjee}\end{center}

In this section we prove Theorem~\ref{thm:extension}, that is used to show that $\Gamma_U$ is a section to $\Psi$ over $U$ in Proposition~\ref{prop:BEsection}.

\begin{setup}\label{setup:ribbon}
Let $\Ct$ be a ribbon of arithmetic genus $g$ and conormal bundle $L_C$ over a smooth curve $C$ of genus $g_C$. Then $\Ct$ is defined by an extension
\[
0\to L_C\to V\to K_C\to0
\]
where $V$ is a rank two vector bundle on $C$. Further $V=\Omega_{\Ct}|_C$. Now the data of the extension canonically gives
\begin{enumerate}
\item a ruled surface $Y=\PP(V)\to C$ and
\item an embedding $C\hookrightarrow Y=\PP(V)$ as a section of the ruled surface $Y$ defined by the surjection $V\to K_C$.
\end{enumerate}
Note that $\deg L_C=2g_C-1-g$, since $\chi(\OO_{\Ct})=\chi(\OO_C)+\chi(L_C)$. 
\end{setup}

%\begin{definition}\label{def:invariant}
%Every ruled surface on a fixed curve $C$ can be uniquely expressed as $Y=\PP(V_n)$ where $V_n$ is a normalized vector bundle, i.e.\ $h^0(V_n)\neq0$ and $h^0(V_n\otimes A^{-1})=0$ for every line bundle $A$ of positive degree. The invariant of a ruled surface is defined as $e=-\deg(V_n)$. Let $C_0$ denote the numerical class of $\OO_{\PP(V_n)}(1)$ and $f$ denote the numerical class of any fiber.
%\end{definition}

\begin{lemma}\label{lem:OSC}
In Set-up~\ref{setup:ribbon}, the line bundle $\OO_Y(C)=\OO_{\PP(V)}(1)\otimes p^*L_C^{-1}$.

\end{lemma}

\begin{proof}
Note that $\mathcal{O}_{\mathbb{P}(V)}(1)|_C = K_C$. Further since $C$ is a section, we have that $\mathcal{O}_Y(C) = \mathcal{O}_{\mathbb{P}(V)}(1) \otimes p^*M$ for some line bundle $M$ on $C$. From the exact sequence
$$0 \to \mathcal{O}_{\mathbb{P}(V)}(1) \otimes \mathcal{O}_{\mathbb{P}(V)}(-C) \to \mathcal{O}_{\mathbb{P}(V)}(1) \to \mathcal{O}_{\mathbb{P}(V)}(1)|_C = K_C \to 0$$
we get 
$$0 \to p^*M^{-1} \to \mathcal{O}_{\mathbb{P}(V)}(1) \to \mathcal{O}_{\mathbb{P}(V)}(1)|_C = K_C \to 0$$

Pushing forward, we have 

$$0 \to M^{-1} \to V \to K_C \to 0$$

This gives, $M^{-1} = L_C$. Hence $\mathcal{O}_Y(C) = \mathcal{O}_{\mathbb{P}(V)}(1) \otimes p^*L_C^{-1}$.

%The section $C$ is the zero locus of the composite $p^*L_C\to p^*V\to\OO_{\PP(V)}(1)$, which vanishes exactly at the points where the quotient of $V$ is $V\to K_C$; hence $\OO_Y(C)=\OO_{\PP(V)}(1)\otimes p^*L_C^{-1}$. We have that $V_n=V\otimes A^{-1}$ where $A$ is the maximum degree line subbundle of $V$. Hence $\OO_{\PP(V)}(1)=\OO_{\PP(V_n)}(1)\otimes p^*A$, and $\OO_Y(C)=\OO_{\PP(V)}(1)\otimes p^*L_C^{-1}=\OO_{\PP(V_n)}(1)\otimes p^*(A\otimes L_C^{-1})$. Now $\deg(A\otimes L_C^{-1})=\deg(A)-\deg(L_C)=2g_C-2-b(\Ct)+g-(2g_C-1)=g-b(\Ct)-1=\frac12(g+e-1)$.

\end{proof}

\begin{theorem}\label{thm:extension}
In Set-up~\ref{setup:ribbon}, assume that $\Ct$ is non-split. For $\yt\in\Ext^1(\Omega_Y,K_Y)$ let $\Yt_{\yt}$ be the $K$-trivial carpet on $Y$ of class $\yt$.
\begin{enumerate}
\item There is a unique $\yt_0\in\Ext^1(\Omega_Y,K_Y)$, $\widetilde y_0 \neq 0$, such that the embedding $i\colon C\hookrightarrow Y$ extends to a closed embedding $\Ct\hookrightarrow\Yt_{\yt_0}$ inducing the identity $K_Y|_C=L_C$ on conormal bundles.
%\item Suppose $\OO_Y(C)$ is very ample, let $Y\hookrightarrow\PP^M$ be the embedding by the complete linear series of $\OO_Y(C)$, $M+1=h^0(\OO_Y(C))$, let $H\subset\PP^M$ be the hyperplane with $Y\cap H=C$, and let $\Yt\hookrightarrow\PP^M$ be the unique $K$-trivial carpet extending the embedding $Y\hookrightarrow\PP^M$ \cite[Proposition~3.2]{DM25}. Then $\Yt\cap H$ is a ribbon on $C$ with conormal bundle $L_C$, and $\Yt\cap H\cong\Ct$ as ribbons on $C$.
\item Let $g_C=0$ and let $\Ct\subset\PP^{M-1}$ be a canonically embedded ribbon, $M=g$, so that $C\subset\PP^{M-1}$ is a rational normal curve. Embed $Y$ in $\PM$ by the complete linear series of $\OO_Y(C)$ and identify $\PP^{M-1}$ with a hyperplane $H\subset\PM$ so that $Y\cap H=C$. If $\Yt\subset\PM$ is the K3 carpet on $Y$, then
\[
\Yt\cap H=\Ct
\]
in $\PP^{M-1}=H$.
\end{enumerate}
\end{theorem}

\begin{proof}
\emph{Proof of (1).} Let $C$ be a smooth curve and $L_C$ a line bundle on $C$. Fix a ribbon $\Ct$ on $C$ with conormal bundle $L_C$. Then $\Ct$ defines an element $[\Ct]\in\Ext^1(\Omega_C,L_C)$. Set $V=\Omega_{\Ct}|_C$. Then we have the exact sequence
\[
0\to L_C\to V\to\Omega_C\to0 .
\]
Let $Y=\PP V$. The exact sequence above gives an embedding
\[
i\colon C\hookrightarrow Y
\]
which is a section of the projection $p\colon Y\to C$. We have
\[
K_Y=\OO_{\PP(V)}(-2)\otimes p^*(K_C\otimes\det V)
\]
and hence
\[
K_Y|_C=\Omega_C^{-2}\otimes(\Omega_C\otimes L_C)\otimes K_C=L_C .
\]
A $K$-trivial carpet $\Yt$ on $Y$ is given by an element $\yt\in\Ext^1(\Omega_Y,K_Y)$. To show the existence of a $K$-trivial carpet as claimed in part (1), we want to choose $\Yt=\Yt_{\yt}$ so that the embedding
\[
C\xhookrightarrow{\ i\ }Y
\]
extends to
\[
\Ct\xhookrightarrow{\ \tilde i\ }\Yt
\]
inducing the identity on conormal bundles. This is equivalent to the equality $\bar i[\Ct]=r(\yt)$ in the following diagram:
\begin{equation}\label{diag:ext}
\begin{tikzcd}
& \Ext^1(\Omega_Y,K_Y)\arrow[d,"r"]\\
\Ext^1(\Omega_C,L_C)\arrow[r,"\bar i"] & \Ext^1(\Omega_Y|_C,L_C)
\end{tikzcd}
\end{equation}
where $r$ is the map induced by the restriction of the short exact sequence $0\to K_Y\to\Omega_{\Yt}|_Y\to\Omega_Y\to0$ to $C$ and $\bar i\colon\Ext^1(\Omega_C,L_C)\to\Ext^1(\Omega_Y|_C,L_C)$ is the map induced by the embedding $i\colon C\hookrightarrow Y$.

Now the projection $p\colon Y\to C$ induces a map $p^*\Omega_C\to\Omega_Y$ which when restricted to the section $i\colon C\hookrightarrow Y$ induces a map $\Omega_C\to\Omega_Y|_C$. This further induces the map $\bar p\colon\Ext^1(\Omega_Y|_C,L_C)\to\Ext^1(\Omega_C,L_C)$ which is a left inverse to the map $\bar i$: $\bar p\circ\bar i=\id$. Thus the conormal sequence $0\to N_{C/Y}^{\vee}\to\Omega_Y|_C\to\Omega_C\to0$ is split by $\Omega_C\to\Omega_Y|_C$, and
\begin{equation}\label{eq:extsplit}
\Ext^1(\Omega_Y|_C,L_C)=\Ext^1(N_{C/Y}^{\vee},L_C)\oplus\Ext^1(\Omega_C,L_C),
\end{equation}
where $\bar i$ is the inclusion of the second summand, $\bar p$ the projection onto it, and the projection onto the first summand is the map
\[
r_N\colon\Ext^1(\Omega_Y|_C,L_C)\to\Ext^1(N_{C/Y}^{\vee},L_C)
\]
induced by $N_{C/Y}^{\vee}\hookrightarrow\Omega_Y|_C$. So it suffices to prove that there is a unique $\yt_0\in\Ext^1(\Omega_Y,K_Y)$ satisfying the two conditions
\begin{equation}\label{eq:twoconditions}
\text{(a)}\ \ [\Ct]=\bar p\,r(\yt_0),\qquad\qquad\text{(b)}\ \ r_N\,r(\yt_0)=0 .
\end{equation}

\emph{Condition (a).} We have to describe the $\yt$ for which $[\Ct]$ is the image of $\yt$ under the map
\[
\Ext^1(\Omega_Y,K_Y)\xrightarrow{\ r\ }\Ext^1(\Omega_Y|_C,L_C)\xrightarrow{\ \bar p\ }\Ext^1(\Omega_C,L_C)
\]

We have a diagram as follows:
\begin{equation}\label{diag:B}
\begin{tikzcd}
\Ext^1(\Omega_{Y/C},K_Y)\arrow[r]\arrow[d] & \Ext^1(\Omega_Y,K_Y)\arrow[r]\arrow[d,"r"] & \Ext^1(p^*\Omega_C,K_Y)=k\arrow[d]\\
\Ext^1(\Omega_{Y/C}|_C,L_C)\arrow[r,hook] & \Ext^1(\Omega_Y|_C,L_C)\arrow[r,"\bar p"] & \Ext^1(\Omega_C,L_C)
\end{tikzcd}
\end{equation}
whose vertical maps are restriction to $C$.
Hence the map the map
\[
\Ext^1(\Omega_Y,K_Y)\xrightarrow{\ r\ }\Ext^1(\Omega_Y|_C,L_C)\xrightarrow{\ \bar p\ }\Ext^1(\Omega_C,L_C)
\]
also factors as
\[
\Ext^1(\Omega_Y,K_Y)\to\Ext^1(p^*\Omega_C,K_Y)\to\Ext^1(\Omega_C,L_C).
\]
We first show that the map $\Ext^1(\Omega_Y,K_Y)\to\Ext^1(p^*\Omega_C,K_Y)$ is surjective. Considering the exact sequence
\[
0\to p^*\Omega_C\to\Omega_Y\to\Omega_{Y/C}\to0 ,
\]
it is enough to show that $\Ext^2(\Omega_{Y/C},K_Y)=0$. By Serre duality,
\[
\Ext^2(\Omega_{Y/C},K_Y)=\Hom(\OO_Y,\Omega_{Y/C})^{\vee}=0
\]
since $\Omega_{Y/C}$ is $\OO_{\PP(V)}(-2)\otimes p^*(\det V)$. First note that since $p_*K_Y=0$, we have
\[
\Ext^1(p^*\Omega_C,K_Y)=H^1(K_Y\otimes p^*T_C)=H^0(R^1p_*K_Y\otimes T_C)=H^0(\OO_C),
\]
because
\[
R^1p_*K_Y=R^1p_*(K_{Y/C}\otimes p^*K_C)=K_C .
\]
So $\Ext^1(p^*\Omega_C,K_Y)=k$. We claim that the image of $1\in\Ext^1(p^*\Omega_C,K_Y)=H^1(\Omega_{Y/C})$ (note that $K_Y\otimes p^*T_C=\Omega_{Y/C}$) in $\Ext^1(\Omega_C,L_C)=H^1(T_C\otimes L_C)$ is $[\Ct]$. This follows from the following general claim.

\begin{claim}\label{claim:euler}
Let $0\to L\to V\to A\to0$ be an extension where $L$ and $A$ are line bundles on $C$. Consider the section $i\colon C\hookrightarrow\PP(V)$ defined by the extension. This induces a map $\Omega_{\PP(V)/C}\to\Omega_{\PP(V)/C}|_C$. Then the image of $1\in H^1(\Omega_{\PP(V)/C})=k$ inside $H^1(\Omega_{\PP(V)/C}|_C)=H^1(L\otimes A^{-1})$ is the class of the extension.
\end{claim}

\emph{Proof of the claim.} $1\in H^1(\Omega_{\PP V/C})$ is represented by the Euler sequence
\[
0\to\Omega_{\PP V/C}\to p^*V\otimes\OO_{\PP V}(-1)\to\OO_{\PP V}\to0 .
\]
Restricted to $C$ via $i$, it is
\[
0\to L\otimes A^{-1}\to V\otimes A^{-1}\to\OO_C\to0 ,
\]
which is the class of the extension twisted by $A^{-1}$. This proves the claim.

For $\yt\in\Ext^1(\Omega_Y,K_Y)$ let $\yt_h\in\Ext^1(p^*\Omega_C,K_Y)=k$ be its image. By the claim,
\begin{equation}\label{eq:pbar-r}
\bar p\,r(\yt)=\yt_h\,[\Ct] ,
\end{equation}
and since $\Ct$ is non-split, $[\Ct]\ne0$; so condition (a) holds for $\yt$ if and only if $\yt_h=1$. The kernel of the map $\Ext^1(\Omega_Y,K_Y)\to\Ext^1(p^*\Omega_C,K_Y)$ is
\[
\Ext^1(\Omega_{Y/C},K_Y)=H^1(T_{Y/C}\otimes K_Y)=H^1(p^*K_C)=H^1(K_C)=k .
\]
Hence the set of $\yt$ satisfying condition (a) is the affine line
\[
B:=\{\yt\in\Ext^1(\Omega_Y,K_Y):\ \yt_h=1\},
\]
the inverse image of $1\in\Ext^1(p^*\Omega_C,K_Y)$.

\emph{Condition (b).} Consider again the diagram:
\begin{equation}
\begin{tikzcd}
\Ext^1(\Omega_{Y/C},K_Y)\arrow[r]\arrow[d] & \Ext^1(\Omega_Y,K_Y)\arrow[r]\arrow[d,"r"] & \Ext^1(p^*\Omega_C,K_Y)=k\arrow[d]\\
\Ext^1(\Omega_{Y/C}|_C,L_C)\arrow[r,hook] & \Ext^1(\Omega_Y|_C,L_C)\arrow[r,"\bar p"] & \Ext^1(\Omega_C,L_C)
\end{tikzcd}
\end{equation}
whose vertical maps are restriction to $C$. Let $B\subset\Ext^1(\Omega_Y,K_Y)$ be the affine line above, and let $K$ be the kernel of $r_N\colon\Ext^1(\Omega_Y|_C,L_C)\to\Ext^1(N_{C/Y}^{\vee},L_C)=k$. We show that $r(B)\cap K$ consists of a single point, i.e.\ that $r_N\circ r$ vanishes at exactly one point of $B$.

Two distinct elements of $B$ differ by a nonzero element of $\Ext^1(\Omega_{Y/C},K_Y)$. This is because
\[
\Hom(p^*\Omega_C,K_Y)=H^0(K_Y\otimes p^*\Omega_C^{*})=0 .
\]
Now note that the conormal sequence
\[
0\to N_{C/Y}^{\vee}\to\Omega_Y|_C\to\Omega_C\to0
\]
has a splitting given by
\[
0\to p^*\Omega_C|_C=\Omega_C\to\Omega_Y|_C\to\Omega_{Y/C}|_C\to0 ,
\]
the restriction to $C$ of $0\to p^*\Omega_C\to\Omega_Y\to\Omega_{Y/C}\to0$. Therefore $\Omega_{Y/C}|_C\cong N_{C/Y}^{\vee}$, and $N_{C/Y}^{\vee}\hookrightarrow\Omega_Y|_C\to\Omega_{Y/C}|_C$ is the identity $N_{C/Y}^{\vee}\to N_{C/Y}^{\vee}$. Hence, for $w\in\Ext^1(\Omega_{Y/C},K_Y)$, $r_N\,r(w)$ is the image of $w$ under the left vertical map of \eqref{diag:B},
\[
\Ext^1(\Omega_{Y/C},K_Y)\to\Ext^1(\Omega_{Y/C}|_C,L_C)=\Ext^1(N_{C/Y}^{\vee},L_C),
\]
and it is enough to show that this map is injective. Tensoring
\[
0\to\OO_Y(-C)\to\OO_Y\to\OO_C\to0
\]
with $\Omega_{Y/C}^{*}\otimes K_Y$ we get
\[
0\to\Omega_{Y/C}^{*}\otimes K_Y\otimes\OO_Y(-C)\to\Omega_{Y/C}^{*}\otimes K_Y\to\Omega_{Y/C}^{*}|_C\otimes L_C\to0 .
\]
Since $\Omega_{Y/C}=K_Y\otimes p^*K_C^{*}$, we have $\Omega_{Y/C}^{*}=K_Y^{*}\otimes p^*K_C$, and hence, by Lemma~\ref{lem:OSC},
\[
\Omega_{Y/C}^{*}\otimes K_Y\otimes\OO_Y(-C)=p^*K_C\otimes\OO_Y(-C)=\OO_{\PP V}(-1)\otimes p^*(K_C\otimes L_C),
\]
whose $H^1$ vanishes because $R^jp_*\OO_{\PP V}(-1)=0$ for all $j$. Taking cohomology of the last exact sequence gives the injectivity of $H^1(\Omega_{Y/C}^{*}\otimes K_Y)\to H^1(\Omega_{Y/C}^{*}|_C\otimes L_C)$, i.e.\ of the left vertical map of \eqref{diag:B}.

Therefore $r_N\circ r$ is injective on $B$. Being an injective affine map from the line $B$ to the line $\Ext^1(N_{C/Y}^{\vee},L_C)=k$, it is bijective, so exactly one point $\yt_0$ of $B$ satisfies condition (b). This proves (1).

%\emph{Proof of (2).} Note that since by Lemma~\ref{lem:OSC}, part (1), $\OO_Y(C)=\OO_{\PP(V)}(1)\otimes p^*L_C^{-1}$, we have the existence of a unique $K$-trivial carpet $\Yt\hookrightarrow\PP^M$ extending the embedding $Y\hookrightarrow\PP^M$, and it is non-split \cite[Proposition~3.2]{DM25}; let $\yt\ne0$ be its class. Let $h$ be a linear form defining $H$, so that $h|_Y$ is the section of $\OO_Y(C)$ with zero scheme $C$. Since $h$ is a nonzerodivisor on $\OO_Y$ and on $K_Y$, it is a nonzerodivisor on $\OO_{\Yt}$, and $K_Y\cap h\OO_{\Yt}=hK_Y$. Hence the kernel of $\OO_{\Yt\cap H}=\OO_{\Yt}/h\OO_{\Yt}\to\OO_Y/h\OO_Y=\OO_C$ is $K_Y/hK_Y=K_Y\otimes\OO_C=L_C$, and it has square zero. Thus $\Ct':=\Yt\cap H$ is a ribbon on $C$ with conormal bundle $L_C$, $\Ct'\cap Y=C$, and the map $K_Y|_C\to I_{C/\Ct'}=L_C$ induced by $\Ct'\subset\Yt$ is the identity. As in the proof of (1), $r(\yt)=\bar i[\Ct']$; hence, by \eqref{eq:extsplit} and \eqref{eq:pbar-r}, $r_N\,r(\yt)=0$ and $[\Ct']=\yt_h[\Ct]$. By the proof of condition (b), $r_N\circ r$ is injective on $\Ext^1(\Omega_{Y/C},K_Y)=\ker(\yt\mapsto\yt_h)$; as $\yt\ne0$ and $r_N\,r(\yt)=0$, we get $\yt_h\ne0$. Therefore $[\Ct']$ and $[\Ct]$ are proportional and nonzero, and $\Ct'\cong\Ct$ as ribbons on $C$.

\emph{Proof of (2).} Since $\Ct$ is non-split, $V\otimes L_C^{-1}\cong\OO_{\PP^1}(a)\oplus\OO_{\PP^1}(a')$ with $a,a'\ge1$ (an extension $0\to\OO_C\to V\otimes L_C^{-1}\to K_C\otimes L_C^{-1}\to0$ with $a=0$ is split), so $\OO_Y(C)=\OO_{\PP(V\otimes L_C^{-1})}(1)$ is very ample; and $h^0(\OO_Y(C))=M+1$, since $H^1(\OO_Y)=0$ and $N_{C/Y}=K_C\otimes L_C^{-1}=\OO_{\PP^1}(g-1)$.

Let $\Yt_{\yt_0}$ and $\Ct\subset\Yt_{\yt_0}$ be as in (1). Then $\Ct$ is a Cartier divisor on $\Yt_{\yt_0}$ with $\OO_{\Yt_{\yt_0}}(\Ct)|_Y=\OO_Y(C)$, and from
\[
0\to K_Y(C)\to\OO_{\Yt_{\yt_0}}(\Ct)\to\OO_Y(C)\to0
\]
and $H^i(K_Y(C))=H^{2-i}(\OO_Y(-C))^{\vee}=0$ for all $i$ (Lemma~\ref{lem:OSC}: $\OO_Y(-C)=\OO_{\PP V}(-1)\otimes p^*L_C$) we get $H^0(\OO_{\Yt_{\yt_0}}(\Ct))\cong H^0(\OO_Y(C))$. Let
\[
\varphi\colon\Yt_{\yt_0}\to\PP^M_1:=\PP\big(H^0(\OO_{\Yt_{\yt_0}}(\Ct))^{\vee}\big)
\]
be the morphism defined by the complete linear series of $\OO_{\Yt_{\yt_0}}(\Ct)$. Its restriction to $Y$ is the embedding by $|\OO_Y(C)|$, and $\varphi^{-1}(H_1)=\Ct$ for the hyperplane $H_1$ defined by the section of $\OO_{\Yt_{\yt_0}}(\Ct)$ vanishing on $\Ct$.

We claim that $\varphi$ is a closed embedding. It is finite (proper and injective on points), so it is enough to show that $\OO_{\PP^M_1}\to\varphi_*\OO_{\Yt_{\yt_0}}$ is surjective. Its image is a subsheaf of algebras mapping onto $\OO_Y=\varphi_*\OO_{\Yt_{\yt_0}}/K_Y$, so it is everything if and only if it contains $K_Y$, i.e.\ if and only if the map $I_Y\to K_Y$ induced by $\OO_{\PP^M_1}\to\varphi_*\OO_{\Yt_{\yt_0}}$ is surjective. Since $K_Y^2=0$ this map factors through $N_Y^{\vee}=I_Y/I_Y^2$, so it is a section of $N_Y\otimes K_Y$. But $H^0(N_Y\otimes K_Y)=k\,s$ with $s$ nowhere vanishing by \cite[Proposition~1.7]{GP97}. Thus the map is $\lambda s$ for some $\lambda\in k$, and it is surjective if $\lambda\ne0$. If $\lambda=0$, then $I_Y$ maps to $0$, so $\OO_{\PP^M_1}\to\varphi_*\OO_{\Yt_{\yt_0}}$ factors through $\OO_Y$, i.e.\ $\varphi$ factors through a retraction of $Y\subset\Yt_{\yt_0}$. But a retraction splits $0\to K_Y\to\OO_{\Yt_{\yt_0}}\to\OO_Y\to0$ as $\OO_Y$-algebras, so $\Yt_{\yt_0}$ would be the split carpet, contradicting $\widetilde y_0 \neq 0$. Hence $\varphi$ is a closed embedding; $\varphi(Y)\subset\PP^M_1$ is embedded by $|\OO_Y(C)|$, $\varphi(\Yt_{\yt_0})$ is a K3 carpet on it, and $\varphi(\Yt_{\yt_0})\cap H_1=\varphi(\Ct)$.

It is clear that $\varphi(\Ct)\subset H_1=\PP^{g-1}$ is canonically embedded.

Now the given $\Ct\subset\PP^{M-1}$ and $\varphi(\Ct)\subset H_1$ are two embeddings of $\Ct$ by the complete linear series of $\omega_{\Ct}$, so there is a projectivity $T\colon H_1\to\PP^{M-1}$ with $T\circ\varphi|_{\Ct}$ the given inclusion; extend $T$ to a projectivity $\hat T\colon\PP^M_1\to\PM$ with $\hat T(H_1)=H$. Then $\hat T\circ\varphi|_Y$ and the given embedding $Y\subset\PM$ are two embeddings of $Y$ by the complete linear series $|\OO_Y(C)|$ which agree on $C$, so they differ by a projectivity $\gamma$ of $\PM$ with $\gamma|_C=\id$, hence $\gamma|_H=\id$ as $C$ spans $H$. Finally $\gamma\hat T\varphi(\Yt_{\yt_0})$ is a K3 carpet on $Y\subset\PM$, hence is $\Yt$ \cite[Proposition~1.7]{GP97}, and
\[
\Yt\cap H=\gamma\hat T\big(\varphi(\Yt_{\yt_0})\cap H_1\big)=\gamma\hat T\varphi(\Ct)=\gamma(\Ct)=\Ct .
\]
\end{proof}

%\section{The compatibility \texorpdfstring{$h\circ f=2g$}{h o f = 2g}}\label{app:hfg}

\end{document}